\documentclass[smallextended,referee,envcountsect]{svjour3}
\smartqed
\usepackage[breaklinks, colorlinks=true]{hyperref}
\usepackage{amsmath,amssymb}
\usepackage{graphicx}
\usepackage{parskip}
\usepackage{mathtools}
\usepackage[margin=2cm, bottom=2.4cm]{geometry}
\usepackage{academicons, xcolor}
\usepackage{float}
\usepackage{orcidlink}
\usepackage{soul}
\usepackage[normalem]{ulem}
\usepackage{cancel}

\newcommand{\A}{\mathcal A}
\renewcommand{\H}{\mathcal H}
\newcommand{\R}{\mathbb R}
\newcommand{\Fix}{\mathrm{Fix}}
\newcommand{\Zer}{\text{Zer}}

\newcommand{\Proj}{\text{Proj}}

\newcommand{\dist}{\text{dist}}

\newcommand{\limty}[1]{\lim_{#1\to\infty}}

\newcommand{\qqbox}[1]{\qquad\hbox{#1}\qquad}

\newtheorem{hypothesis}{Hypothesis}[section]

\journalname{JOTA}

\begin{document}

\title{Krasnoselskii-Mann Iterations: Inertia, Perturbations and Approximation}

\author{Daniel Cortild \orcidlink{0000-0002-3278-1716} \and Juan Peypouquet \orcidlink{0000-0002-8551-0522}}

\date{Submitted: January 31, 2024, Accepted: 11 December 2024, Published online: 22 January 2025}

\institute{Daniel Cortild \at
    University of Groningen \\
    d.cortild@rug.nl
\and
    Juan Peypouquet \at
    University of Groningen \\
    j.g.peypouquet@rug.nl
}

\maketitle

\begin{abstract}
This paper is concerned with the study of a family of fixed point iterations combining relaxation with different inertial (acceleration) principles. We provide a systematic, unified and insightful analysis of the hypotheses that ensure their weak, strong and linear convergence, either matching or improving previous results obtained by analysing particular cases separately. We also show that these methods are robust with respect to different kinds of perturbations--which may come from computational errors, intentional deviations, as well as regularisation or approximation schemes--under surprisingly weak assumptions. Although we mostly focus on theoretical aspects, numerical illustrations in image inpainting and electricity production markets reveal possible trends in the behaviour of these types of methods.


\end{abstract} 


\section{Introduction}

Let $\H$ be a real Hilbert space. Krasnoselskii-Mann iterations \cite{Krasnoselskii,Mann} approximate fixed points of a (quasi) nonexpansive operator $T:\H\to \H$, by means of the update rule
$$x_{k+1}=(1-\lambda_k)x_k+\lambda_kTx_k,$$
where $\lambda_k\in(0,1)$ is a {\it relaxation} parameter. They were independently introduced by Mann in 1953 \cite{Mann}, with $\lambda_k=\frac{1}{k+1}$, and by Krasnoselskii in 1955 \cite{Krasnoselskii}, with $\lambda_k\equiv \frac{1}{2}$. Their weak convergence was established in \cite[Krollar 2.1]{Schaefer1957berDM} for any constant parameter $\lambda_k\equiv\lambda\in(0,1)$, and then in \cite[Corollary 3]{GROETSCH1972369} for variable relaxation parameters satisfying $\sum\lambda_k(1-\lambda_k)=\infty$. Krasnoselskii-Mann iterations are central to numerical optimization and variational analysis, where many problems can be reduced to finding fixed points of appropriate operators. Many known splitting optimization algorithms are special instances.

On the other hand, the consideration of physical principles has proven to be a useful technique in optimization. The concept of {\it momentum} was first introduced by Polyak in 1964 \cite{POLYAK19641}, who showed that the Heavy Ball method accelerates convergence in certain problems. Although originally proposed for gradient descent methods, it may be extended to Krasnoselskii-Mann iterations \cite{alvarez2001inertial,AccKM,dong2022new}, giving
\begin{equation*}
\begin{cases}
y_k&=\quad x_k+\alpha_k(x_k-x_{k-1}) \\
x_{k+1}&=\quad (1-\lambda_k)y_k+\lambda_kTx_k,
\end{cases}
\end{equation*}
where $\alpha_k$ is an {\it acceleration} parameter. The idea of momentum was later reinterpreted by Nesterov in 1983 \cite{Nesterov}, also to accelerate the convergence of gradient methods. Since then, many algorithms have been improved by the addition of this more popular acceleration step, especially in the context of convex optimization \cite{beck2009fast,doi:10.1137/0802032}, although this is not our emphasis here. Nesterov's acceleration scheme has also been used in fixed-point theory and variational analysis \cite{boct2015inertial,Dong,iyiola2021new,mainge2008convergence,FierroMaulenPeypouquet,MOUDAFI2003447,shehu2018convergence}, under the form
\begin{equation*}
\begin{cases}
z_k&=\quad x_k+\beta_k(x_k-x_{k-1}) \\
x_{k+1}&=\quad (1-\lambda_k)z_k+\lambda_kTz_k,
\end{cases}
\end{equation*}
where $\beta_k$ is an {\it inertial} parameter (observe the differences and similarities with Polyak's approach).  These two interpretations of inertia can be combined into a more general algorithm \cite{Dong2022,Dong2018,dong2021general,gebregiorgis2023convergence}
\begin{equation*}
\begin{cases}
y_k&=\quad x_k+\alpha_k(x_k-x_{k-1}) \\
z_k&=\quad x_k+\beta_k(x_k-x_{k-1}) \\
x_{k+1}&=\quad (1-\lambda_k)y_k+\lambda_kTz_k.
\end{cases}
\end{equation*}
This not only provides a unified setting for the study of the classical inertial methods described above, but its versatility also suggests new ones. For instance, for $\alpha_k\equiv 0$, we obtain the {\it reflected} Krasnoselskii-Mann iterations \cite{Dong2022,Dong2018,doi:10.1080/10556788.2021.1924715,moudafi2018reflected}, inspired by the reflected gradient method \cite{malitsky2015projected}.

The purpose of this work is threefold:
\begin{enumerate}
    \item First, to provide a systematic and unified analysis of the hypotheses that ensure the convergence of the sequences produced by means of inertial Krasnoselskii-Mann iterations. In doing so, we either match or extend the range of admissible parameters known to date, which had previously been obtained by analysing different particular cases separately.
    \item Next, to establish the extent to which these iterations are stable with respect to perturbations, which could be due to computational or approximation errors, or to deviations purposely introduced in order to enhance different aspects of the algorithms' approximation power (further commentary in Subsection \ref{sec:errors}). 
    \item Finally, to account for {\it diagonal} algorithms \cite{peypouquet2009asymptotic,attouch2011prox,attouch2011coupling,noun2013forward,peypouquet2012coupling} represented by a sequence of operators, a situation that typically arises when the iterative procedure is coupled with regularisation or approximation strategies.
\end{enumerate}

To this end, we study the behaviour of sequences generated iteratively by the set of rules
\begin{equation}\label{alg:0}
\begin{cases}
y_k&=\quad x_k+\alpha_k(x_k-x_{k-1}) +\varepsilon_k \\
z_k&=\quad x_k+\beta_k(x_k-x_{k-1}) +\rho_k \\
x_{k+1}&=\quad (1-\lambda_k)y_k+\lambda_kT_kz_k +\theta_k,
\end{cases}
\end{equation}
where $T_k:\H\to\H$, $\alpha_k,\beta_k,\lambda_k\in[0,1]$, and $\varepsilon_k,\rho_k, \theta_k\in\H$ for $k\ge 1$, and $x_0,x_1\in\H$ (which we set equal, for simplicity). 

The paper is organized as follows: Section \ref{sec:convergence} contains the convergence  results, where weak, strong and linear convergence is established. In Section \ref{S:comments}, we provide more insight into the hypotheses concerning the relationships between the parameters (Subsection \ref{SS:parameters}), how sharp our conditions are (Subsection \ref{ssec:tight}), as well as a more detailed comparison with other results found in the literature, which we improve both in the exact and the perturbed cases (Subsection \ref{SS:connection}), before addressing the extension to families of operators not sharing a fixed point, which accounts for approximation or regularization procedures (Subsection \ref{sec:general}) and a few comments of the nature of the perturbation sequences (Subsection \ref{sec:errors}). Although our work is mostly concerned with theoretical aspects of these methods, we include some numerical examples in Section \ref{S:numerics} that illustrate how the different kinds of inertial schemes behave in an image inpainting problem (Subsection \ref{sec:tos}), and the search for a Nash-Cournot oligopolistic equilibrium in electricity production markets (Subsection \ref{sec:eq}). Finally, Appendix \ref{sec:postponed} contains the postponed technical proofs, and some auxiliary results are gathered in Appendix \ref{sec:aux}.

\section{Convergence Results} \label{sec:convergence}

Throughout this paper, $\H$ is a real Hilbert space with inner product $\langle\cdot, \cdot\rangle$ and induced norm $\|\cdot\|$. Strong (norm) and weak convergence of sequences will be denoted by $\to$ and $\rightharpoonup$, respectively. The set of fixed points of an operator $T\colon \H\to \H$ is $\Fix(T)\coloneqq \{x\in \H\colon Tx=x\}$. 

\subsection{Weak Convergence} \label{SS:weak_convergence}

We recall an operator $T$ is \textit{quasi-nonexpansive} if $\Fix(T)\neq \emptyset$, and $|Tx-p|\le |x-p|$ for all $x\in \H$ and $p\in \Fix(T)$.

In order to simplify the statements of the results below, let us summarize our standing assumptions on the parameter sequences, where we define $\mu_k \coloneqq  (1-\lambda_k)\alpha_k+\lambda_k\beta_k$:

\begin{hypothesis} \label{H:weak}
The sequences $(\alpha_k)$, $(\beta_k)$, $(\lambda_k)$ and $(\mu_k)$, along with the constants $\alpha=\inf\alpha_k$, $A=\sup\alpha_k$, $\lambda=\inf\lambda_k$, $\Lambda=\sup\lambda_k$ and $M=\sup\mu_k$, satisfy: $\alpha,A, M\in [0, 1)$, $\lambda,\Lambda\in(0,1)$, $0\le\beta_k\le 1$,  and $\mu_k\le\mu_{k+1}$ for all $k\ge 1$. 
\end{hypothesis}

\begin{remark}
The sequence $(\mu_k)$ is chosen by the user via the parameters of the algorithm. For the Heavy Ball momentum ($\beta_k\equiv 0$) we have $\mu_k=(1-\lambda_k)\alpha_k$; for Nesterov's acceleration ($\alpha_k=\beta_k$), it reduces to $\mu_k=\alpha_k$; and for the reflected acceleration ($\alpha_k\equiv 0$), we have $\mu_k=\lambda_k\beta_k$. The technical hypothesis that $\mu_k\le\mu_{k+1}$ for all $k\ge 1$ is trivially satisfied in the constant case, which is the most common in practice, except for the {\it Fast Iterative Shrinkage-Thresholding Algorithm} (FISTA, \cite{Nesterov,beck2009fast}), which also assumes $\mu_k=\alpha_k\le\alpha_{k+1}=\mu_{k+1}$, for all $k\ge 1$. 
\end{remark}

The following results establish asymptotic properties of the sequences generated by Algorithm \eqref{alg:0}, when the parameter sequences satisfy Hypothesis \ref{H:weak} and the compatibility condition
\begin{equation} \label{eq:tbr}
\sup_{k\ge 1}\left[\frac{}{}\!(1-\lambda_k)\alpha_k(1+\alpha_k)+\lambda_k\beta_k(1+\beta_k) +\nu_k\alpha_k(1-\alpha_k)-\nu_{k-1}(1-\alpha_{k-1})\right] < 0,
\end{equation}
where $\nu_k \coloneqq \lambda_k^{-1}-1$.

\begin{remark} \label{R:constraints}
Inequality \eqref{eq:tbr} cannot hold if 1 is a limit point of either $(\alpha_k)$ or $(\lambda_k)$. This is the motivation to assume that $\sup\alpha_k<1$ and $\sup\lambda_k<1$ in Hypothesis \ref{H:weak}. Also, since we are interested in asymptotic results, the supremum in Inequality \eqref{eq:tbr} can be replaced by an upper limit. Further insight into this inequality is provided in Subsection \ref{SS:parameters}.
\end{remark}

\begin{proposition}\label{prop:2}
Let $T_k\colon \H\to \H$ be a family of quasi-nonexpansive operators such that $F\coloneqq \bigcap_{k\ge 1}\Fix(T_k)\neq \emptyset$, let Hypothesis \ref{H:weak} and Inequality \eqref{eq:tbr} hold, and assume the error sequences $(\varepsilon_k)$, $(\rho_k)$ and $(\theta_k)$ belong to $\ell^1(\H)$. If $(x_k,y_k,z_k)$ is generated by Algorithm \eqref{alg:0}, then
\begin{equation} \label{E:summability_residuals}
\sum_{k=1}^\infty \|x_k-x_{k-1}\|^2<\infty\qqbox{and} \sum_{k=1}^\infty \|T_kz_k-y_k\|^2<\infty.
\end{equation}
In particular, $\limty{k}\|x_k-x_{k-1}\|=\limty{k}\|T_kz_k-y_k\|=0.$ Moreover, $\limty{k}\|x_k-p\|$ exists for all $p\in F$.
\end{proposition}

The proof is deferred to Appendix \ref{sec:prop2}.

\begin{remark}
The boundedness of the sequence $(x_k)$, as well as the square-summability of the residuals \eqref{E:summability_residuals} are part of the conclusion of Proposition \ref{prop:2}, not hypotheses. This contrasts with several recent works found in the literature, even in the unperturbed case (see Subsection \ref{SS:connection}).
\end{remark}

\begin{remark}
In the context of Proposition \ref{prop:2}, $\min_{i=1,\dots,k}\|T_iz_i-y_i\|^2=o(\frac{1}{k})$.   
\end{remark}

A family $(T_k)$ of operators is \textit{asymptotically demiclosed} (at $0$) if, for every sequence $(u_k)$ in $\H$, such that $u_k\rightharpoonup u$ and $T_ku_k-u_k\to 0$, it follows that $u\in \bigcap_{k\ge 1}\Fix(T_k)$.

\begin{theorem}\label{thm:weak}
Let $T_k\colon \H\to \H$ be a family of quasi-nonexpansive operators such that $F\coloneqq \bigcap_{k\ge 1}\Fix(T_k)\neq \emptyset$, let Hypothesis \ref{H:weak} and Inequality \eqref{eq:tbr} hold, and suppose the error sequences $(\varepsilon_k)$, $(\rho_k)$ and $(\theta_k)$ belong to $\ell^1(\H)$. Assume, moreover, that $(I-T_k)$ is asymptotically demiclosed at $0$. If $(x_k,y_k,z_k)$ is generated by Algorithm \eqref{alg:0}, then $(x_k,y_k,z_k)$ converges weakly to $(p^*,p^*,p^*)$, with $p^*\in F$.
\end{theorem}

\begin{proof}
Since $x_{k+1}-x_k\to 0$ (Proposition \ref{prop:2}), $\varepsilon_k\to 0$ and $\rho_k\to 0$, the definition of $y_k$ and $z_k$, implies that $(x_k)$, $(y_k)$ and $(z_k)$ have the same set of weak limit points. Moreover, since $(\alpha_k)$ and $(\beta_k)$ are bounded,
\begin{equation*}
y_k-z_k=(\alpha_k-\beta_k)(x_k-x_{k-1})+\varepsilon_k-\rho_k\to 0.
\end{equation*}
Also, $\lim_{k\to \infty}\|x_k-p\|$ exists for all $p\in F$, and 
\begin{equation*}
(I-T_k)z_k=(y_k-T_kz_k)-(y_k-z_k)\to 0.
\end{equation*} 
The asymptotic demiclosedness of $(I-T_k)$ then implies that every weak limit point $(z_k)$ must belong to $F$, and the same is true for every weak limit point of $(x_k)$. By Opial's Lemma \cite{Opial}, $(x_k)$--as well as $(y_k)$ and $(z_k)$--must converge weakly to some $p^*\in F$. \qed 
\end{proof}

The connection with specific instances of this method, along with a comparison with the results found in the literature are discussed in Subsection \ref{SS:connection}.

\subsection{Strong Convergence}\label{sec:strong}

We now analyze the case where each $T_k$ is $q_k$-quasi-contractive, which means that $\|T_kx-p\|\le q_k\|x-p\|$ for all $x\in \H$ and $p\in \Fix(T_k)$, and some $q_k<1$.

As before we define $\nu_k=\lambda_k^{-1}-1$, and we additionally define $Q_k=1-\lambda_k+\lambda_kq_k^2$. We are now in a position to show the strong convergence of the sequences generated by Algorithm \eqref{alg:0}, when the parameter sequences satisfy
\begin{equation}\label{eq:tbr_strong}
\sup_{k\ge 1}\left[\frac{}{}\!(1-\lambda_k)\alpha_k(1+\alpha_k)+\lambda_kq_k^2\beta_k(1+\beta_k)+\nu_k\alpha_k(1-\alpha_k)-Q_k\nu_{k-1}(1-\alpha_{k-1})\right]<0,
\end{equation}
which is reduced to \eqref{eq:tbr} when $q_k\equiv 1$ (which corresponds to the quasi-nonexpansive case). Notice that Remark \ref{R:constraints} remains pertinent.

\begin{theorem}\label{thm:strong}
Let $T_k\colon \H\to \H$ be a family of $q_k$-quasi-contractive operators with $q_k\le q<1$, and such that $\Fix(T_k)\equiv\{p^*\}$. Let Hypothesis \ref{H:weak} and Inequality \eqref{eq:tbr_strong} hold, and assume the error sequences $(\varepsilon_k)$, $(\rho_k)$ and $(\theta_k)$ belong to $\ell^2(\H)$. 
If $(x_k,y_k,z_k)$ is generated by Algorithm \eqref{alg:0}, then 
$(x_k,y_k,z_k)$ converges strongly to $(p^*,p^*,p^*)$. Moreover, $\sum_{k=1}^\infty\|x_k-p^*\|^2<\infty$, and $\min_{j=1,\dots,k}\|x_j-p^*\|^2=o(1/k)$, as $k\to\infty$. If, in addition, $\varepsilon_k\equiv \rho_k\equiv \theta_k\equiv 0$, then
$$\|x_k-p^*\|^2\le Q^k\frac{\|x_1-p^*\|^2}{(1-\Lambda)(1-A)}$$
for all $k\ge 1$, where $Q:=\sup_{k\ge 0}Q_k$.
\end{theorem}

The proof is deferred to Appendix \ref{sec:thmstrong}.
\section{Discussion and Implications} \label{S:comments}

\subsection{The Relationships Between the Parameters} \label{SS:parameters}

To simplify the exposition and fix the ideas, let us restrict ourselves to the case of constant parameters, namely $\alpha_k\equiv\alpha$, $\beta_k\equiv\beta$ and $\lambda_k\equiv\lambda$. Inequality \eqref{eq:tbr} is reduced to 
$$\lambda(1-\lambda)\alpha(1+\alpha)+\lambda^2\beta(1+\beta) +(1-\lambda)\alpha(1-\alpha)-(1-\lambda)(1-\alpha)<0.$$
In other words,
$$(\beta-\alpha)(1+\alpha+\beta)\lambda^2+(1-\alpha+2\alpha^2)\lambda-(1-\alpha)^2<0.$$

For $\beta=0$ (Heavy Ball), this gives
$(1-\lambda)\big(\alpha(1+\alpha)\lambda-(1-\alpha)^2\big)<0.$
Since $\lambda\le 1$, this means that
$$\lambda<\lambda_{HB}(\alpha):=\frac{(1-\alpha)^2}{\alpha(1+\alpha)}.$$
The right-hand side is greater than 1 if $\alpha<1/3$, so there is no constraint on $\lambda$ in those cases!

For $\alpha=\beta$ (Nesterov), the coefficient in the second-order term disappears, and we are left with
$$\lambda<\lambda_N(\alpha):=\frac{(1-\alpha)^2}{1-\alpha+2\alpha^2},$$
which coincides with the hypothesis used in \cite{FierroMaulenPeypouquet} in this particular case. The right hand side decreases from 1 to 0 when $\alpha$ goes from 0 to 1.

For $\alpha=0$ (Reflected), we get
$(1+\beta\lambda)\big((1+\beta)\lambda-1\big)<0,$
which is nothing more than
$$\lambda<\lambda_R(\beta):=\frac{1}{1+\beta},$$
which decreases from 1 to $1/2$ as $\beta$ goes from 0 to 1.

\begin{figure}[H]
    \centering
    \includegraphics[width=0.6\linewidth]{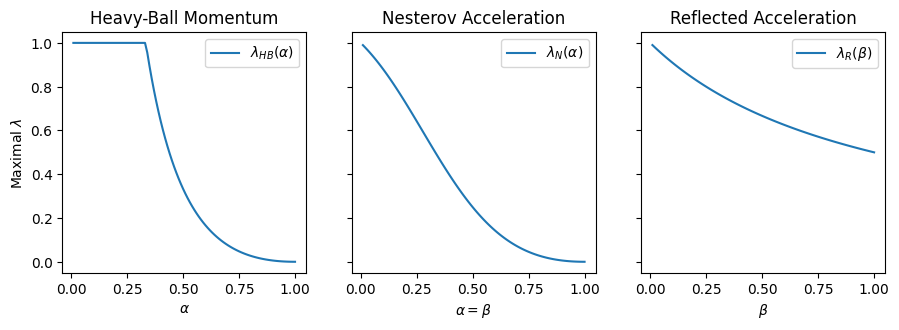}
    \caption{Upper bounds on $\lambda$ for Heavy-Ball (left), Nesterov (center) and Reflected (right) acceleration. }
    \label{fig:lamb}
\end{figure}

The functions $\lambda_{HB}$, $\lambda_N$ and $\lambda_R$ are depicted in Figure \ref{fig:lamb}.

Extending by continuity, the constraint can be written explicitly in the general case as
\begin{equation} \label{eq:lambda_alpha_beta}
\lambda<\lambda(\alpha,\beta):=
\frac{\sqrt{(1-3\alpha)^2+4\beta(1+\beta)(1-\alpha)^2}-1+\alpha-2\alpha^2}{2(\beta-\alpha)(1+\alpha+\beta)}.
\end{equation}

\begin{figure}[H]
    \centering
    \includegraphics[width=0.4\linewidth]{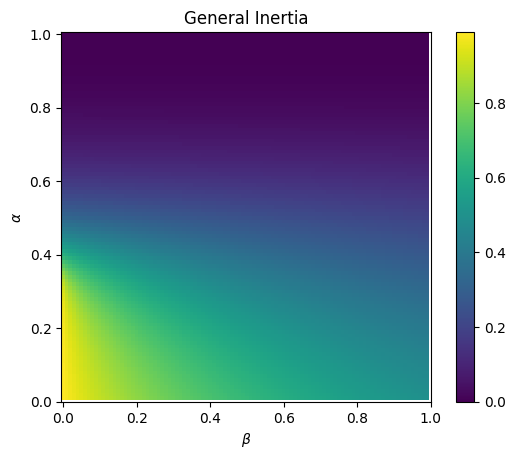}
    \caption{Upper bounds on $\lambda$ for general inertia as a function of $(\alpha, \beta)$.}
    \label{fig:general3}
\end{figure}

The function $(\alpha,\beta)\mapsto\lambda(\alpha,\beta)$ is shown in Figure \ref{fig:general3}.

\subsection{Tightness of the Constraints on the Parameters}\label{ssec:tight}

\begin{figure}[H]
    \centering
    \includegraphics[width=0.5\linewidth]{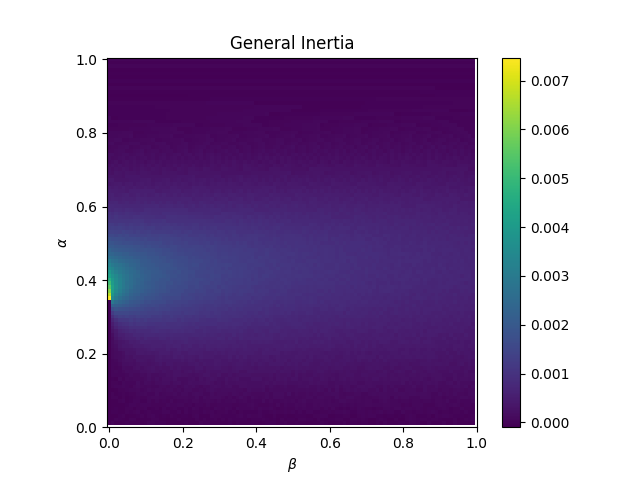}
    \caption{Values of $\tilde\lambda(\alpha,\beta)- \lambda(\alpha,\beta)$ for $(\alpha,\beta)\in[0,1)\times [0, 1]$.}
    \label{fig:enter-label}
\end{figure}

In order to assess the sharpness of the condition given by Inequality \eqref{eq:tbr}, we analyse the convergence of the iterations defined by Algorithm \eqref{alg:0} for a family of nonexpansive operators with a simple structure. Consider $\lambda_k\equiv\lambda\in(0,1)$, $\alpha_k\equiv\alpha\in[0,1)$, $\beta_k\equiv\beta\in[0,1)$ and $T_k\equiv T_\phi$, where $T_\phi:\R^2\to\R^2$ is a counterclockwise rotation on the plane by an angle of $\phi\in(0,\pi]$. To simplify the computations, we identify $\R^2$ with $\mathbb C$ in order to represent $T_\phi$, using complex multiplication, as $z\mapsto e^{i\phi}z$. The algorithm is then expressed as 
$$x_{k+1}=\left[(1-\lambda)(1+\alpha)-\lambda(1+\beta)e^{i\phi}\right]x_k-\left[(1-\lambda)\alpha+\lambda \beta e^{i\phi}\right]x_{k-1}=\omega x_k + \delta x_{k-1},$$
where we have written $\delta=(1-\lambda)\alpha+\lambda\beta e^{i\phi}$ and $\omega=(1-\lambda)(1+\alpha)-\lambda(1+\beta)e^{i\phi}$. Passing to the product space, we can write the iterations in matrix form as
$$\left(\begin{array}{c}x_{k+1} \\ x_k\end{array}\right)=\left(\begin{array}{cc}\omega & -\delta \\ 1 & 0\end{array}\right) \left(\begin{array}{c}x_{k} \\ x_{k-1}\end{array}\right).$$
The algorithm converges if, and only if, both complex eigenvalues of the matrix have modulus less than 1. For each $\alpha,\beta,\phi$, let 
$\tilde \lambda(\alpha, \beta, \phi)$ be the supremum over all $\lambda>0$ for which the algorithm converges. Then, set $$\tilde \lambda(\alpha, \beta)\coloneqq \inf_{\phi}\tilde \lambda (\alpha, \beta, \phi),$$ 
so that
$$\lambda(\alpha,\beta)\le \tilde\lambda(\alpha,\beta)\le \tilde \lambda (\alpha, \beta, \phi)$$
for every $\alpha,\beta,\phi$. If $\tilde\lambda(\alpha,\beta)-\lambda(\alpha,\beta)=0$ for all $\alpha,\beta$, then Inequality \eqref{eq:tbr} is tight. To avoid the tedious algebraic work of computing an analytic expression for $\tilde\lambda(\alpha,\beta)$, we parameterize the family $T_\phi$ over $\phi\in \{j\pi/30\colon j=1, \ldots, 30\}$. Table \ref{tab:tightness} shows different measures of the difference between $\lambda$ and $\tilde\lambda$, namely
$$\|\tilde\lambda-\lambda\|_1=\int_{\mathcal D}\big|\tilde\lambda(\alpha,\beta)-\lambda(\alpha,\beta)\big|\,d(\alpha,\beta)\quad\hbox{and}\quad \|\tilde\lambda-\lambda\|_\infty={\sup}_{(\alpha,\beta)\in\mathcal D}\big|\tilde\lambda(\alpha,\beta)-\lambda(\alpha,\beta)\big|,$$
for the corresponding domain $\mathcal D$. Figure \ref{fig:enter-label} shows the point-wise values of $\tilde\lambda-\lambda$ over $[0, 1)\times [0, 1]$. The very small values may be partially due to computer precision, although the highest values coincide with the {\it kink} shown in Figure \ref{fig:lamb} in our bound for the Heavy Ball.

\begin{table}[H]
    \centering
    \begin{tabular}{|c|c|c|c|} \hline
        & $\mathcal D$ & $\|\tilde\lambda-\lambda\|_1$ & $\|\tilde\lambda-\lambda\|_\infty$ \\ \hline
        Heavy-Ball & $\{(\alpha, 0)\colon \alpha\in [0, 1)\}$   
                & 0.00086   & 0.00747 \\ \hline
        Nesterov   & $\{(\alpha, \alpha)\colon \alpha\in [0, 1)\}$      
                & 0.00047   & 0.00122 \\ \hline
        Reflected  & $\{(0, \beta)\colon \beta\in [0, 1]\}$   
                & 0.00005   & 0.00010 \\ \hline
        General    & $\{(\alpha, \beta)\colon \alpha\in [0, 1), \beta\in [0, 1]\}$  
                & 0.00001   & 0.00747  \\ \hline
    \end{tabular}
    \caption{Distance to tightness of Inequality \eqref{eq:tbr} or, more particularly, Inequality \eqref{eq:lambda_alpha_beta}, measured in terms of $L^1$ and $L^\infty$ norms.}
    \label{tab:tightness}
\end{table}

\subsection{Our Results in Perspective} \label{SS:connection}

Convergence of inertial methods has been studied by a number of researchers in different, and mostly less general contexts. We discuss here how the results established above relate to previously known ones.

\subsubsection{Exact Methods}
Several special cases of Algorithm \eqref{alg:0} in the unperturbed case $\varepsilon_k\equiv\rho_k\equiv \theta_k\equiv 0$ have already been studied, namely:

$\bullet$ Nesterov's acceleration corresponds to  $\alpha_k=\beta_k$. Weak convergence has been established in \cite{lorenz2015inertial} for the forward-backward method, assuming square-summability of the residuals. This hypothesis was proved unnecessary in \cite{iutzeler2019generic}, and convergence was proved under hypotheses equivalent to ours in the constant case. Also, \cite[Theorem 5]{boct2015inertial}, \cite[Theorem 3.1]{shehu2018convergence} and \cite[Proposition 3.1 and Theorem 3.1]{Dong}, show weak convergence
under hypotheses on the parameters that are similar, but more involved. Indeed, the author of \cite{Dong} remarks that {\it these
conditions are too complicated to determine an upper bound for the inertial sequence in a simple
way, even if the coefficient $\lambda_k$ is known. Moreover, in the case of $\lambda_k\equiv 0.5$, they are
undesirably restrictive}\footnote{For readability, the notation in this quote has been adapted to match ours, and some misprints have been corrected.}. The said hypotheses were simplified in \cite[Theorem 2.1]{iyiola2021new}, but are still more restrictive than ours, since they constrain $\alpha\in\left[0,\frac{1}{3}\right)$. All of the above consider $T_k\equiv T$, with $T$ nonexpansive. The case of a family of operators was studied in \cite[Theorems 3.1 and 3.2]{mainge2008convergence}, both in the nonexpansive and in the firmly quasi-nonexpansive cases, but also assuming a summability condition on the residuals, which we have proved, based on the arguments in \cite{FierroMaulenPeypouquet}, to be unnecessary, even under perturbations. On a different note, Nesterov's acceleration has also been added to algorithms of extragradient type \cite[Theorems 3.1 and 3.2]{thong2018inertial}.

$\bullet$ The case $\beta_k\equiv 0$ (Heavy Ball momentum), was studied in \cite{alvarez2001inertial,AccKM}, where convergence is obtained assuming square-summability and boundedness of the residuals, respectively, both impractical hypotheses. This was solved in \cite[Theorem 1]{dong2022new}, where convergence is proved under assumptions similar to ours.

$\bullet$ Reflected acceleration, which corresponds to $\alpha_k\equiv 0$, was studied in \cite[Theorem 4]{Dong2018} and \cite[Theorem 5.4]{Dong2022}, under slightly stronger assumptions on the parameters. The very particular case $\beta_k\equiv 1$ was analysed in \cite[Theorem 3.1]{doi:10.1080/10556788.2021.1924715}. On the other hand, \cite[Proposition 2.1]{moudafi2018reflected} includes an additional projection step. Weak convergence is obtained for Lipschitz pseudo-contrative mappings, and strong convergence for Lipschitz strongly monotone mappings with monotonicity constant strictly larger than $1$. 

$\bullet$ In the remaining cases (although still considering $\varepsilon_k\equiv\rho_k\equiv \theta_k\equiv 0$ and $T_k\equiv T$), weak convergence of Algorithm \eqref{alg:0} was established in \cite[Theorem 1]{Dong2018} and \cite[Theorem 5.1]{Dong2022}, assuming that $(\lambda_k)$ is constant, both $(\alpha_k)$ and $(\beta_k)$ are nondecreasing, and an additional condition (in line with \cite{boct2015inertial,shehu2018convergence,Dong}), which the authors qualify as {\it complicated} and {\it restrictive}, in \cite{Dong2022}. An online selection of the relaxation parameters is studied in \cite[Theorem 3]{dong2021general}. Much stronger hypotheses (on the parameters, the operator and the residuals) are used in \cite[Theorem 2]{gebregiorgis2023convergence}. An application to three-operator splitting is given in \cite[Theorem 1]{wen2019two}, under {\it ad hoc} assumptions. Finally, a multi-step inertial Krasnoselskii-Mann algorithm is studied in \cite[Theorem 4.1]{dong2019mikm}, but also assuming summability of the residuals.

\subsubsection{Perturbed Algorithms}

For Nesterov's acceleration, convergence under $\ell^1$ perturbations in the Krasnoselskii-Mann step is proved in \cite[Theorem 3.1]{cui2019convergence}. However, they assume the sequence $(x_k)$ to be bounded, which is impractical and, as we show, unnecessary. Similar results were obtained earlier in \cite[Theorems 3.5 and 5.1]{khatibzadeh2015inexact}, assuming both that the generated sequence is bounded and that the residuals are square-summable. An inexact version of FISTA was analyzed in \cite[Theorem 5.1]{attouch2018fast} without any boundedness hypotheses, but assuming that the perturbations satisfy the summability condition $\sum k\|\varepsilon_k\|<+\infty$, which is stronger than ours. Relative error conditions have been accounted for in \cite[Theorem 3.1]{padcharoen2019convergence}, \cite[Theorem 3.6]{alves2020relative} and \cite[Theorem 2.5, 3.3 and 4.4]{alves2020relativeb}, as well as \cite[Theorems 3.2 and 3.3]{ezeora2022inexact}, although the latter is concerned with {\it alternated} inertia. Other approaches include \cite[Algorithm 1]{aujol2015stability}, as well as an inexact multilayer FISTA \cite[Algorithm 2.1]{lauga2023multilevel} and an inexact accelerated forward-backward algorithm \cite[Theorem 3.2]{villa2013accelerated}, where the summability conditions on the perturbations are taylored to the corresponding methods. 

\subsection{Operators Not Sharing a Fixed Point}\label{sec:general}

Our results can also be applied when $\bigcap_{k\ge 1}\Fix(T_k)=\emptyset$, if instead we are interested in the {\it Kuratowski lower limit} of the family $\big(\Fix(T_k)\big)$, which we denote by  $F_\infty$, and consists of all $p_\infty\in\H$ for which there is a sequence $(p_k)$, such that $p_k\in \Fix(T_k)$ for all $k\ge 1$, and $p_k\to p_\infty$. The sequence $(x_k)$ will {\it pursue} these sets, and will finally converge to a point in $F_\infty$. The following simple example in the non-accelerated case can be insightful:

\begin{example}
Let $(p_k)$ be a sequence in $\H$, and let $T_k\equiv p_k$, so that $\Fix(T_k)=\{p_k\}$, for each $k\ge 1$. Define the sequence $(x_k)$ by means of the Krasnoselskii-Mann iterations $x_{k+1}=(1-\lambda)x_k+\lambda T_kx_k=(1-\lambda)x_k+\lambda p_k$, for $k\ge 1$, where $\lambda\in(0,1)$. 
Using Lemma \ref{L:diff}, one can prove: (i) that $p_k$ converges to $p^*$ if, and only if, $x_k$ does; (ii) that $(p_k-p_{k-1})\in\ell^1(\H)$ if, and only if, $(x_k-x_{k-1})\in\ell^1(\H)$ (the polygonal interpolations have finite length); and (iii) that $(p_k-p_{k-1})\in\ell^2(\H)$ if, and only if, $(x_k-x_{k-1})\in\ell^2(\H)$.
\end{example}

We can extend the results of Theorems \ref{thm:weak} and \ref{thm:strong} to this {\it moving set} context. To this end, first fix $p_\infty\in F_\infty$, and set $p_k=\Proj_{\Fix(T_k)}p_\infty$. For each $k$, define $\tilde T_k\colon \H\to \H$ by
\begin{equation*}\label{def:tilde}
\tilde T_k x= T_k(x+p_k-p_\infty)-p_k+p_\infty.
\end{equation*}
Clearly, $p_\infty\in \tilde F\coloneqq \bigcap_{k\ge 1}\Fix(\tilde T_k)$, so
$\Fix(\tilde T_k)\neq\emptyset$, and $\tilde p_k\in\Fix(\tilde T_k)$ if, and only if, $\tilde p_k+p_k-p_\infty\in\Fix(T_k)$. Since $T_k$ is quasi-nonexpansive, for each $x\in \H$, we have
\begin{equation*}
\|\tilde T_kx-\tilde p_k\|=\|T_k(x+p_k-p_\infty)-(\tilde p_k+p_k-p_\infty)\|\le \|(x+p_k-p_\infty)-(\tilde p_k+p_k-p_\infty)\|= \|x-\tilde p_k\|,
\end{equation*}
and so, $\tilde T_k$ is quasi-nonexpansive, as well. Quasi-contractivity is also inherited by $\tilde T_k$ from $T_k$.

By setting $\tilde x_k=x_k-p_{k-1}+p_\infty$, $\tilde y_k=y_k-p_k+p_\infty$, and $\tilde z_k=z_k-p_k+p_\infty$, Algorithm \eqref{alg:0} becomes
\begin{equation*}
\begin{cases}
\tilde y_k &= \quad \tilde x_k+\alpha_k(\tilde x_k-\tilde x_{k-1}) + \tilde \varepsilon_k,\qqbox{where}\tilde\varepsilon_k:=p_{k-1}-p_k+\alpha_k(p_{k-1}-p_{k-2})+\varepsilon_k \\
\tilde z_k &= \quad \tilde x_k+\beta_k(\tilde x_k-\tilde x_{k-1}) + \tilde\rho_k,\qqbox{where}\tilde\rho_k:=p_{k-1}-p_k+\beta_k(p_{k-1}-p_{k-2})+\rho_k \\
\tilde x_{k+1} &= \quad (1-\lambda_k)\tilde y_k+\lambda_k\tilde T_k\tilde z_k+\theta_k.
\end{cases}
\end{equation*}
If $(\varepsilon_k),(\rho_k),(p_k-p_{k-1})\in\ell^r(\H)$, then also $(\tilde \varepsilon_k),(\tilde \rho_k)\in\ell^r(\H)$.

\begin{remark}
It is neither necessary to know $p_\infty$ and the sequence $(p_k)$, nor to construct the operators $\tilde T_k$ and the auxiliary variables $\tilde x_k$, $\tilde y_k$ or $\tilde z_k$. These are merely artifacts to prove convergence.
\end{remark} 

Proposition \ref{prop:2} shows that $\limty{k}\|\tilde x_k-\tilde x_{k-1}\|=\limty{k}\|\tilde T_k\tilde z_k-\tilde y_k\|=0$, and $\limty{k}\|\tilde x_k-p_\infty\|$ exists. As a consequence, $\limty{k}\|x_k-x_{k-1}\|=\limty{k}\|T_kz_k-z_k\|=0$ (because $\limty{k}\|z_k-y_k\|=0$), and $\limty{k}\|x_k-p_\infty\|$ exists. This is true for each $p_\infty\in F_\infty$.
    
Let us say that $(T_k)$ {\it nicely approximates} $F_\infty$ if $u_k\rightharpoonup u$ and $T_ku_k-u_k\to 0$ together imply $u\in F_\infty$. 

\begin{example}
In $\H=\R^N$, this holds if there is a strictly increasing continuous function $\Phi:\R\to\R$ such that $\Phi(0)=0$ and 
$$\dist\big(z,\Fix(T_k)\big)\le \Phi\big(\|z-T_kz\|\big)$$  
for every $k\ge 1$, and $z\in \H$. This is similar to the error bound in \cite{Luo_Tseng}, and can also be understood as the family $(\varphi_k)$ of functions, defined by $\varphi_k(z)=\|z-T_kz\|$, having a common residual function $\Phi$ (see \cite[Section 2.4]{bolte2017error}).    
\end{example}

From Theorem \ref{thm:weak}, we obtain:

\begin{corollary}
Let $T_k\colon \H\to \H$ be a family of quasi-nonexpansive operators that nicely approximates $F_\infty\neq \emptyset$. Let Hypothesis \ref{H:weak} and Inequality \eqref{eq:tbr} hold, and assume the error sequences $(\varepsilon_k)$, $(\rho_k)$, and $(\theta_k)$ belong to $\ell^1(\H)$. Assume moreover that there exists a sequence $(p_k)$ with $p_k\in \Fix(T_k)$ and $(p_k-p_{k-1})\in\ell^1(\H)$. If $(x_k,y_k,z_k)$ is generated by Algorithm \eqref{alg:0}, then $(x_k,y_k,z_k)$ converges weakly to some $(p^*,p^*,p^*)$, with $p^*\in F_\infty$.
\end{corollary}

The preceding discussion also allows us to prove the following extension of Theorem \ref{thm:strong}: 

\begin{corollary}
Let $T_k\colon \H\to \H$ be a family of $q_k$-quasi-contractive operators with $q_k\le q<1$, and such that $\Fix(T_k)=\{p_k\}$, with $p_k\to p^*$ and $(p_k-p_{k-1})\in\ell^2(\H)$. Let Hypothesis \ref{H:weak} and Inequality \eqref{eq:tbr_strong} hold, and assume the error sequences $(\varepsilon_k)$, $(\rho_k)$, and $(\theta_k)$ belong to $\ell^2(\H)$. If $(x_k,y_k,z_k)$ is generated by Algorithm \eqref{alg:0}, then $(x_k,y_k,z_k)$ converges strongly to $(p^*,p^*,p^*)$. Moreover, $\sum_{k=1}^\infty\|x_k-p^*\|^2<\infty$.
\end{corollary}

\subsection{The Role of Perturbations}\label{sec:errors}

A number of different circumstances can give rise to the sequences $(\varepsilon_k)$, $(\rho_k)$ and $(\theta_k)$. In principle, our results do not assume anything about the {\it nature} of the perturbations, as long as they are either in $\ell^1(\H)$ or $\ell^2(\H)$, respectively. The $\ell^1(\H)$ summability is a standard assumption in the analysis of inexact methods in variational analysis; the $\ell^2(\H)$ condition is not nearly as common. Although it is not the purpose of this work to go any deeper into the actual implementation issues concerning these perturbations$-$certainly an important matter$-$, let us briefly comment on a few examples:

\begin{enumerate}
    \item Approximation errors may appear upon the implementation of the fixed-point operator(s). These can have different origins, for instance:
    \begin{enumerate}
        \item Inexactness of numerical methods, when the operators involve integration, operator inversion, or the resolution of (either ordinary or partial) differential equations. In some cases, one can estimate these errors and enforce a prescribed accuracy.
        \item Limited information in data-based models. Although accuracy estimates may exist, the availability of data plays a major role.
    \end{enumerate}
    \item Deviations or bias purposefully introduced in order to enhance different aspects of the algorithm's approximation power, such as:
    \begin{enumerate}
        \item Selection of a particular solution (least norm, closest to initial point, or any other secondary optimization criterion), or in the moving set context of Subsection \ref{sec:general}.
        \item Heuristics that improve performance \cite{sadeghi2023incorporating,sadeghi2021forward,sadeghi2022dwifob}, {\it superiorization} \cite{censor2022superiorization}, and so on.
    \end{enumerate}
    \item Stochastic approximations, either arising from actual uncertainty, or as a means to reduce iteration cost.
    \begin{enumerate}
        \item In the fundamentally uncertain case, a natural course of action to obtain almost sure convergence would be to adapt Lemmas A.3 and A.4 to the context of \cite[Theorem 1]{robbins1971convergence}, as done in \cite{combettes2015stochastic}. We believe that positive results may be obtained in this fashion, and is an interesting direction for future research. 
        \item In the spirit of stochastic gradient descent, the extension of our analysis seems tricky for two reasons: first, the analysis required to deduce {\it convergence in expectation} for several different criteria usually relies on estimations on the expected value of the  objective function, and second, convergence usually requires to {\it neutralise} the effect of the variance, often by recurring to vanishing step sizes (learning rates), which are not allowed with our arguments. In the contracting case, one can overcome the first difficulty, but not necessarily the second one. 
    \end{enumerate}
\end{enumerate}

\section{Numerical Results} \label{S:numerics}

In this section, we illustrate how the different popular acceleration schemes behave in two examples. In Subsection \ref{sec:tos}, we look into the image inpainting problem, using the three operator scheme. Secondly, in Subsection \ref{sec:eq}, we consider a Nash-Cournot oligopolistic equilibrium model in electricity markets.

The experiments are implemented in Python 3.11. In both settings, we shall use an error function $\mathcal R\colon \H\to \R$ defined by
\begin{equation}\label{eq:error}
\mathcal R(X)= \|T(X)-X\|,
\end{equation}
measuring the distance from $X$ to $T(X)$. The convergence results from the previous sections imply that $\mathcal R(X_k)\to 0$ if, and only if, $X_k\to p^*$ where $\Fix(T)=\{p^*\}$. We will terminate the algorithm when the iterates reach a given tolerance $\mathcal R(X_k)\le \varepsilon$, or when a pre-specified maximal number of iterations is reached, after which we consider the algorithm not to have converged.

\subsection{Image Inpainting Problem}\label{sec:tos}

For $A,B\colon \H\rightrightarrows \H$ maximally monotone operators and $C\colon \H\to\H$ a $\tau$-cocoercive operator, we aim to find $\hat x\in \H$ such that 
\begin{equation}\label{problem:zero}
\hat x\in \Zer(A+B+C).
\end{equation}

Such a scheme is called a \textit{three-operator splitting scheme} \cite{DavisYin}, and is equivalent to finding $\bigcap_{k\ge 1}\Fix(T_k)$, where $T_k$ is defined as
\begin{equation*}
T_k\coloneqq I-J_{\rho_k B}+J_{\rho_k A}\circ (2J_{\rho_k B}-I-\rho_k C\circ J_{\rho_k B}).
\end{equation*}
Indeed, $J_{\rho_k B}(\Fix(T_k))= \Zer(A+B+C)$, and the operator $T_k$ is nonexpansive when $(\rho_k) \subset (0,2\tau)$. As such, Algorithm \eqref{alg:0} converges under the given conditions. 

Consider the \textit{image inpainting problem}: We represent an image $X$ of $M$ by $N$ pixels by a tensor in $\mathcal H\coloneqq\R^{M\times N\times 3}$, in which the three layers represent the red, green and blue colour channels. Let $\Omega$ be an element of $\{0,1\}^{M\times N}$ such that $\Omega_{ij}=0$ indicates that the pixel at position $(i,j)$, on all colour channels, has been damaged. Denote by $\mathcal A$ the linear operator that maps an image to an image whose elements in $\Omega$ have been erased. More precisely,
\begin{equation*}
\mathcal A\colon \mathcal H\to \mathcal H, \quad X\mapsto \tilde X,\quad \text{where } \tilde X_{ijk}=\Omega_{ij}\cdot X_{ijk}.
\end{equation*}
The operator $\A$ is a self-adjoint bounded projection map with operator norm $1$. We denote the damaged image by $X_{\text{corrupt}}=\A X$. The objective is to recover an image from $X_{\text{corrupt}}$ that  mainly overlaps on the points where $\Omega_{ij}=1$, and which looks better to the eye, which is obtained by adding the regularization $\|X_{(1)}\|_*+\|X_{(2)}\|_*$, where $X_{(1)}\coloneqq [X_{\cdot \cdot 1}~ X_{\cdot \cdot 2}~ X_{\cdot \cdot 3}]$, $X_{(2)}\coloneqq \left[X_{\cdot \cdot 1}^T~ X_{\cdot \cdot 2}^T~ X_{\cdot \cdot 3}^T\right]^T$ and $\|\cdot\|_*$ denotes the nuclear norm. The image inpainting problem is
\begin{equation*}\label{eq:miniprob}
\min_{X\in \H}\left\{\frac12\|\A X - X_{\text{corrupt}}\|^2 +\sigma \|X_{(1)}\|_*+\sigma \|X_{(2)}\|_*\right\},
\end{equation*}
where $\sigma>0$ is a  \textit{regularisation parameter}.
Since the functions involved are continuous everywhere, this problem can be described by \eqref{problem:zero}, by virtue of the Moreau-Rockafellar Theorem (see, for instance, \cite[Theorem 3.30]{Pey}). To this end, first write $f(X)=\sigma \|X_{(1)}\|_*$, $g(X)=\sigma \|X_{(2)}\|_*$, $h(X)=\frac12 \|X-X_{\text{corrupt}}\|^2$ and $L=\A$. Then, set $A\coloneqq\partial f$, $B\coloneqq\partial g$, both maximally monotone, and $C\coloneqq\nabla (h\circ L)$, which is $\tau/\|L\|_{\text{op}}$-cocoercive. 

The image to be inpainted has dimensions $512\times 512$ pixels. We select a regularisation parameter of $\sigma=0.5$, a tolerance of $\varepsilon = 0.5$, and a maximal number of iterations of $100$. We corrupt the images randomly, with a certain percentage of pixels erased on all colour channels. We always set $X_0=X_1=X_{\text{corrupt}}$.

For simplicity, we set $\rho_k\equiv \rho\in (0,2)$ and $\lambda_k\equiv \lambda\in (0,1)$, and add no perturbations (other than possible rounding errors by the machine). We run multiple versions of the algorithm, corresponding to different inertial schemes: no inertia, Nesterov, Heavy Ball, and reflected. In each case, we pick $\alpha$ and $\beta$ such that Inequality \eqref{eq:tbr} is tight, and then select
\begin{equation*}
    \alpha_k=\left(1-\frac1k\right)\alpha, \quad \beta_k=\left(1-\frac1k\right)\beta.
\end{equation*}

\subsubsection*{The Tests}

\begin{figure}[H]
\begin{center}
\includegraphics[width=1\linewidth]{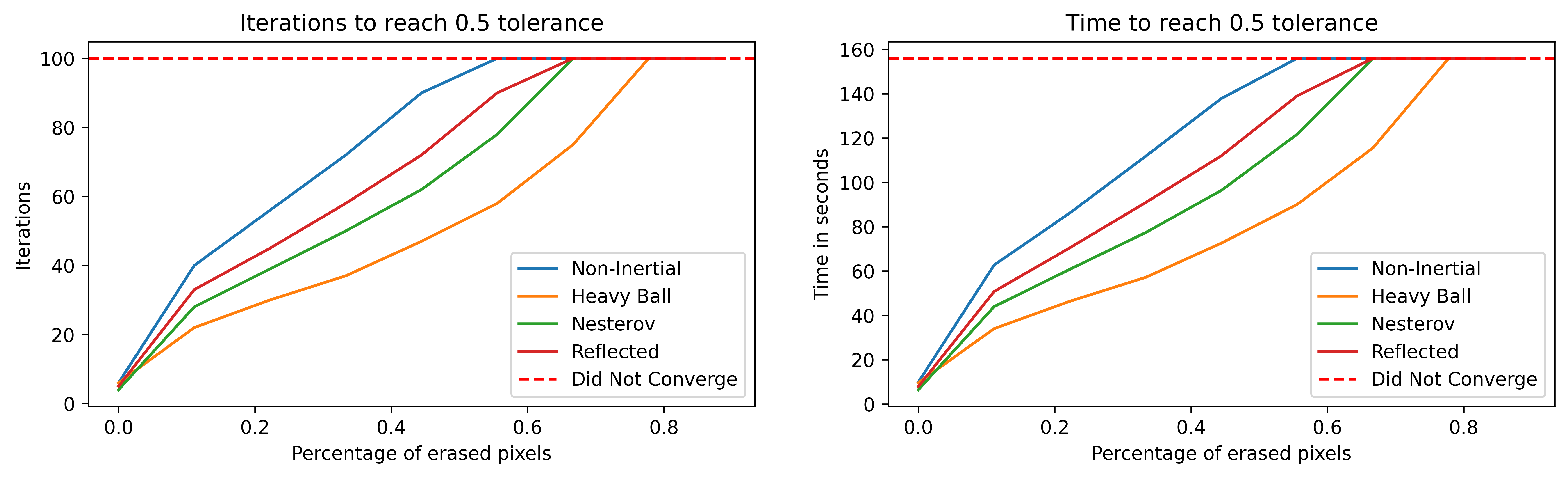}
\end{center}
\caption{Number of iterations and execution time for different ratios of erased pixels, with $\rho=1$ and $\lambda=0.5$.}
\label{fig:2}
\end{figure}

First, we compare the evolution of the number of iterations and the execution time required by our algorithms to inpaint a randomly corrupted image, as a function of the percentage of pixels erased, whilst fixing the step size $\rho=1$ and the relaxation parameter $\lambda=0.5$. The results are shown in Figure \ref{fig:2}. As could be expected, we observe an overall increasing trend. The Heavy Ball acceleration performs overall the best, and is capable of restoring the image with close to $80\%$ of its pixels erased in the given number of iterations, whereas the non-inertial version only is when at most $50\%$ are erased.

\begin{figure}[H]
\begin{center}
\includegraphics[width=1\linewidth]{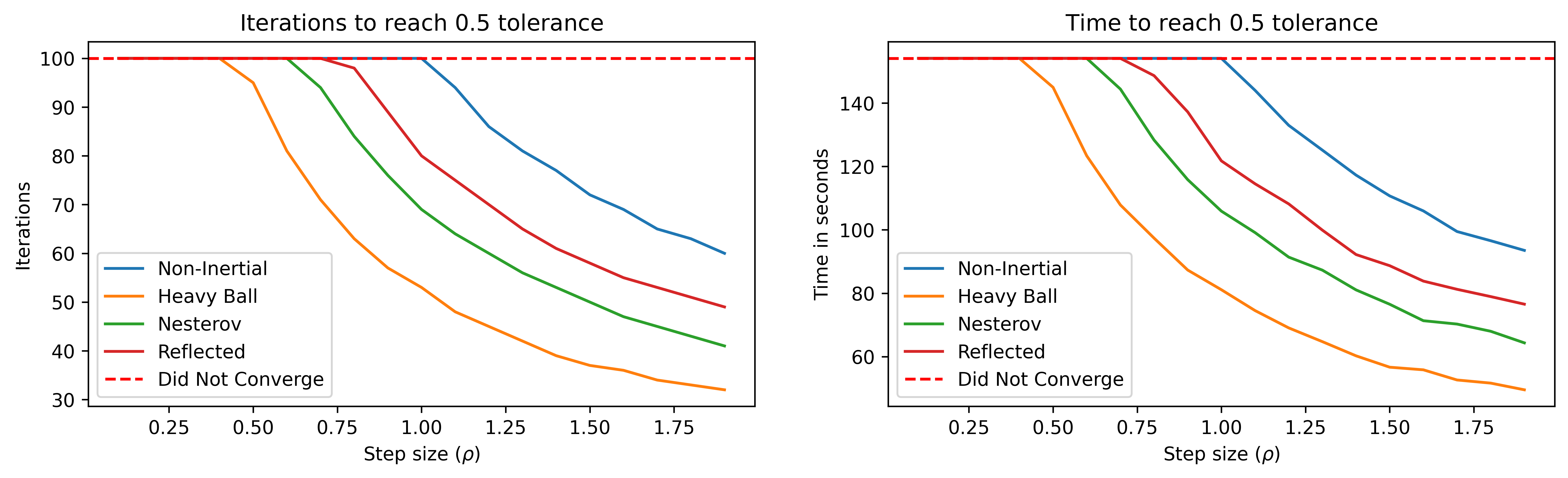}
\end{center}
\caption{Number of iterations and execution time for $\lambda=0.5$ and a ratio of erased pixels of $50\%$, as a function of the step size $\rho$.}
\label{fig:3}
\end{figure}

Next, we fix the relaxation parameter $\lambda=0.5$ and randomly corrupt $50\%$ of the pixels in the image. We iterate over representative values of the step size $\rho\in (0,2)$. The results are shown in Figure \ref{fig:3}. Similar observations as for the percentage with respect to the comparison of the different versions may be made. Additionally, we notice that a larger value of $\rho$ accelerates the convergence.

\begin{figure}[H]
\begin{center}
\includegraphics[width=0.8\linewidth]{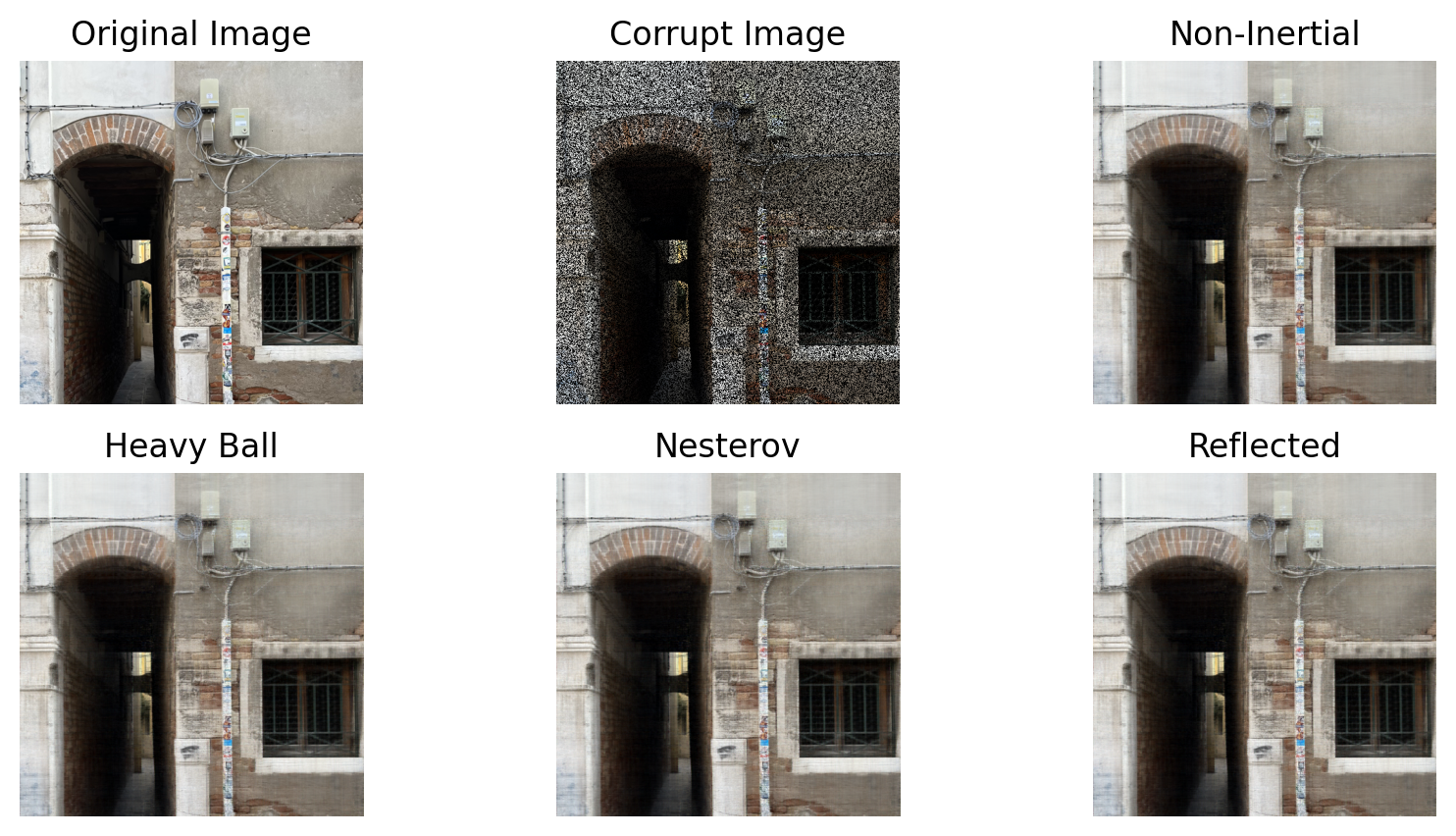}
\end{center}
\caption{Inpainted image with $\lambda=0.8$, $\rho=1.8$, and $\sigma=0.5$.}
\label{fig:5}
\end{figure}

\begin{figure}[H]
    \begin{center}
    \includegraphics[width=0.7\linewidth]{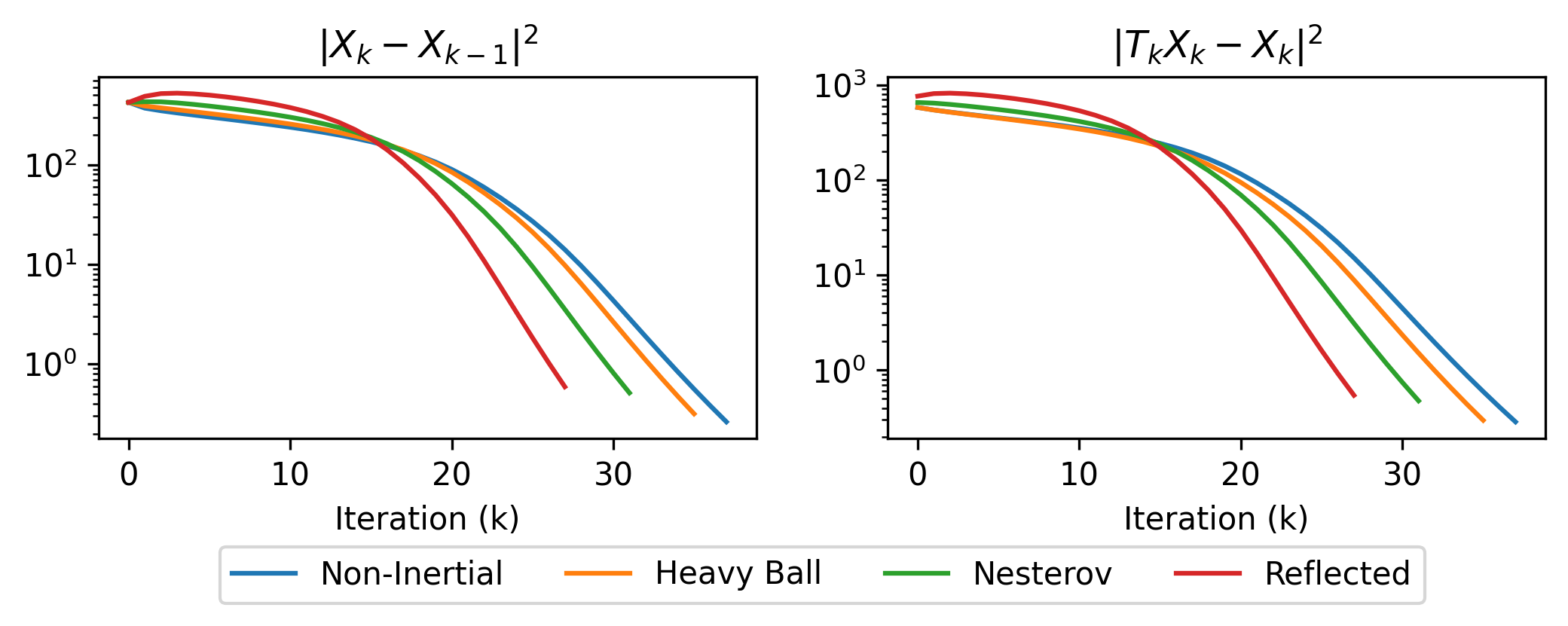}
    \end{center}
    \caption{Convergence plots of the residuals of the iterations.}
    \label{fig:6}
\end{figure}

Finally, visual results for each version of the algorithm are shown in Figure \ref{fig:5}, and the convergence rates in Figure \ref{fig:6}. Here, we considered a corrupted image with $50\%$ of erased pixels, regularisation parameter $\sigma=0.5$, step size $\rho= 1.8$, and relaxation parameter $\lambda= 0.8$.

Despite the similarity between the recovered images in Figure \ref{fig:5}, one can see a quantitative difference in the convergence plots in Figure \ref{fig:6}. The non-inertial version of the algorithm is always outperformed by the three inertial versions, and the reflected acceleration converges faster than the two other inertial variants. We restricted ourselves to these four variants of the algorithm, and leave the interesting$-$and highly challenging!$-$problem of optimising the parameters for future studies.

\subsection{Nash-Cournot Equilibrium Model}\label{sec:eq}

We now consider a Nash-Cournot oligopolistic equilibrium model in electricity production markets \cite{contreras2004numerical}. We consider $m$ companies, and denote by $x_i$ the power generated by company $i$. The generation price for company $i$ is $p_i(s)=\gamma-\beta_i\cdot s$, where $s=\sum_{i=1}^mx_i$. The profit made by company $i$ is $f_i(x)=p_i(s)x_i-c_i(x_i)$, where $c_i(x_i)$ is the cost to generate $x_i$ by company $i$. We denote $C_i$ the strategy set of company $i$ so that $x_i\in C_i$, and define the strategy profile set by $C=C_1\times\cdots\times C_m$.

In Cournot competition, no firms cooperate or collude--each company seeks to maximise their profit and assume the remaining firms do the same. A point $\hat x\in C$ is a \textit{Nash equilibrium} if, for all $i=1, \ldots, m$,  
\[
    f_i(\hat x)=\max_{x_i\in C_i} f_i(\hat x[x_i]),
\]
where $\hat x[x_i]$ represents the vector $\hat x$ whose $i$th component has been replaced by the entry $x_i$. This means that, under strategy profile $\hat x$, no single firm benefits from deviating from the strategy $\hat x$.

We define the Nikaido-Isoda \cite{pjm/1171984836} function $f\colon C\times C\to \R$ by $f(x,y)=\sum_{i=1}^m\left(f_i(x)-f_i(x[y_i])\right)$.
We may now write the above problem as finding $\hat x\in C$ such that $f(\hat x, x)\ge 0$ for all $x\in C$. Assuming the cost functions $c_i$ are convex and differentiable, this can be rewritten as \cite{yen2016algorithm}
\begin{equation}\label{prob:var}
    \langle F(\hat x), x-\hat x\rangle \ge 0\quad \forall x\in C,
\end{equation}
where $F(x)=B x-\Gamma+\nabla \phi(x)$ with
\begin{equation*}
\Gamma\coloneqq (\gamma, \cdots, \gamma)^T, \quad
\tilde B\coloneqq \begin{pmatrix}
    \beta_1 & 0 & \cdots & 0 \\
    0 & \beta_2 & \cdots & 0 \\
    \vdots & \vdots & \ddots & \vdots \\
    0 & 0 & \cdots & \beta_m \\
   \end{pmatrix}, \quad
    B\coloneqq \begin{pmatrix}
         0 & \beta_1 & \cdots & \beta_1 \\
         \beta_2 & 0 & \cdots & \beta_2 \\
         \vdots & \vdots & \ddots & \vdots \\
         \beta_m & \beta_m & \cdots & 0 \\
    \end{pmatrix}, \quad
    \phi(x)\coloneqq x^T\tilde B x+\sum_{i=1}^mc_i(x_i).
\end{equation*}
Variational inequality problems such as Problem \eqref{prob:var} where $C$ is a nonempty, closed and convex set and $F$ is monotone and $L$-Lipschitz have been extensively studied \cite{xu2010iterative,malitsky2015projected,he2018selective,zhao2018iterative,yao2019convergence,dong2021two}. The simplest iterative procedure to solve the above variational inequality is through the projected gradient method, namely
\begin{equation*}\label{eq:VI_Solve}
    x_{k+1}=P_C(x_n-\rho F(x_n)),
\end{equation*}
where $P_C\colon \H\to C$ represents the projection onto $C$. By writing $T_n\equiv T=P_C(I-\rho F)$, we obtain a family of nonexpansive operators whose fixed points represent a solution to Problem \eqref{prob:var}, provided that the step size $\rho$ satisfies $0<\rho<2/L$
\cite{xu2010iterative}. We notice that the method by Malitsky \cite{malitsky2015projected} is a reflected acceleration of the previous, with parameter $\beta_n\equiv1$.

We shall assume the cost functions are quadratic, namely of the form
\begin{equation*}
    c_i(x_i)=\frac12 p_ix_i^2+q_ix_i,
\end{equation*}
where $p_i>0$ such that $L=\|B-p^TI\|$, where $p=(p_1, \ldots, p_m)^T$. The parameters $\beta_i, p_i, q_i$ are selected uniformly at random in the intervals $(0,1]$, $[1,3]$ and $[1,3]$ respectively, and we assume $C_i=[1,40]$, for all $1\le i\le m$. We set $m=8$, $\gamma=200$ and $\rho=1/L$. We select a tolerance of $\varepsilon=10^{-4}$, and a maximum number of iterations of $800$.

As in the previous example, we run multiple versions of the algorithm, corresponding to the different types of inertia, no perturbations and $\alpha_k,\beta_k$ chosen as before.

\subsubsection*{The Tests}

\begin{figure}[H]
\includegraphics[width=\linewidth]{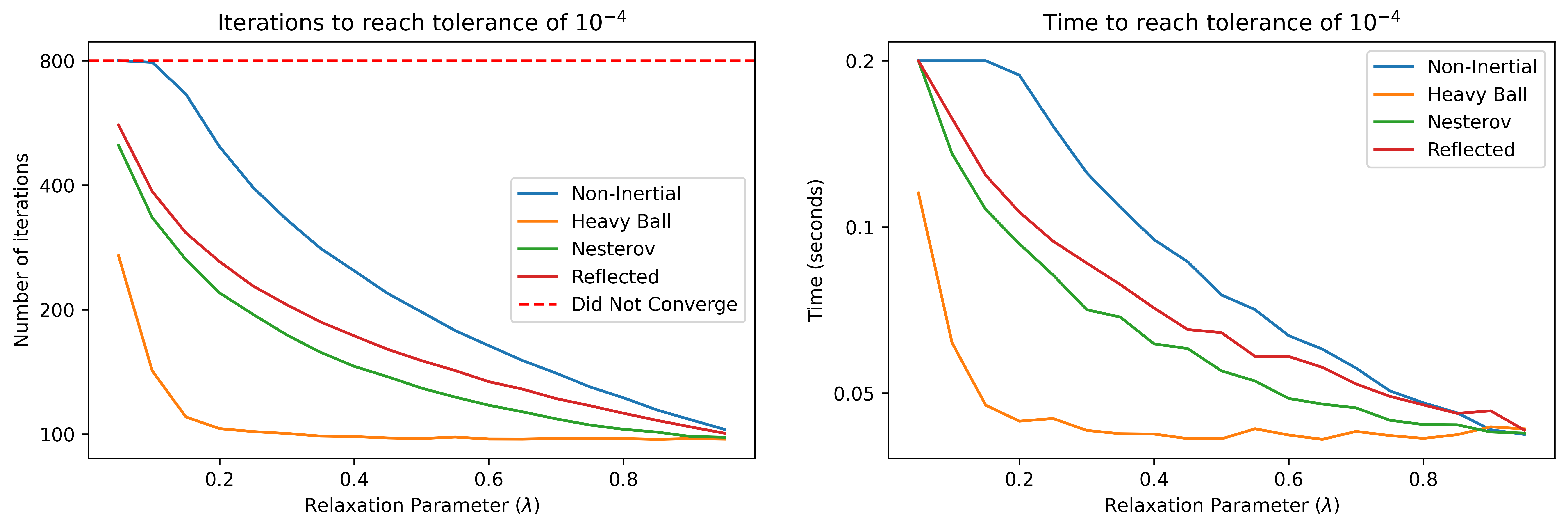}
    \caption{Variation of the relaxation parameter $\lambda$ for $\rho=1/L$.}
    \label{fig:VI_lamb}
\end{figure}

We run the algorithm for various values of the relaxation parameter $\lambda\in(0,1)$. We observe convergence for all the selected values of $\lambda$, and faster convergence for all methods as $\lambda\approx 1$. The Heavy Ball acceleration produces the best results, and the required number of iterations are roughly identical for $\lambda\ge 0.2$. The results are depicted in Figure \ref{fig:VI_lamb}. 

\begin{figure}[H]
    \includegraphics[width=\linewidth]{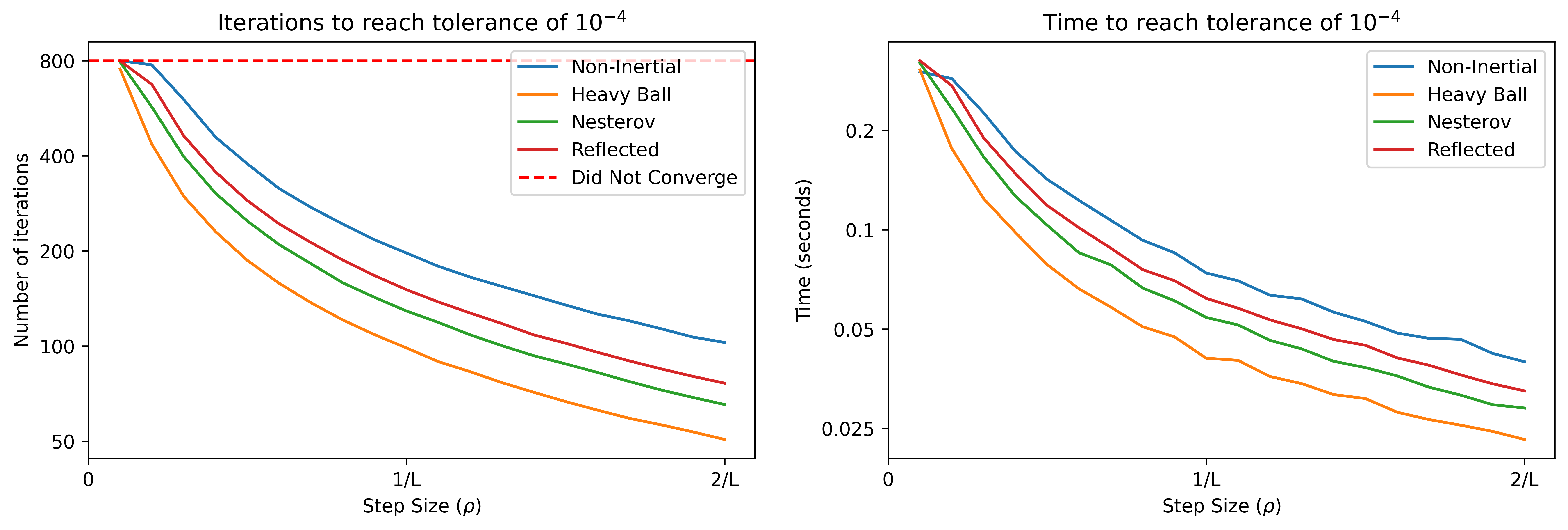}
    \caption{Variation of the step size $\rho\in (0, 2/L)$ for $\lambda=0.5$. }
    \label{fig:VI_rho}
\end{figure}

We fix the relaxation parameter to be $\lambda=0.5$, and execute the algorithm for various step-sizes $\rho\in(0, 2/L)$. The convergence results are plotted in Figure \ref{fig:VI_rho}. 

\begin{figure}[H]
\includegraphics[width=\linewidth]{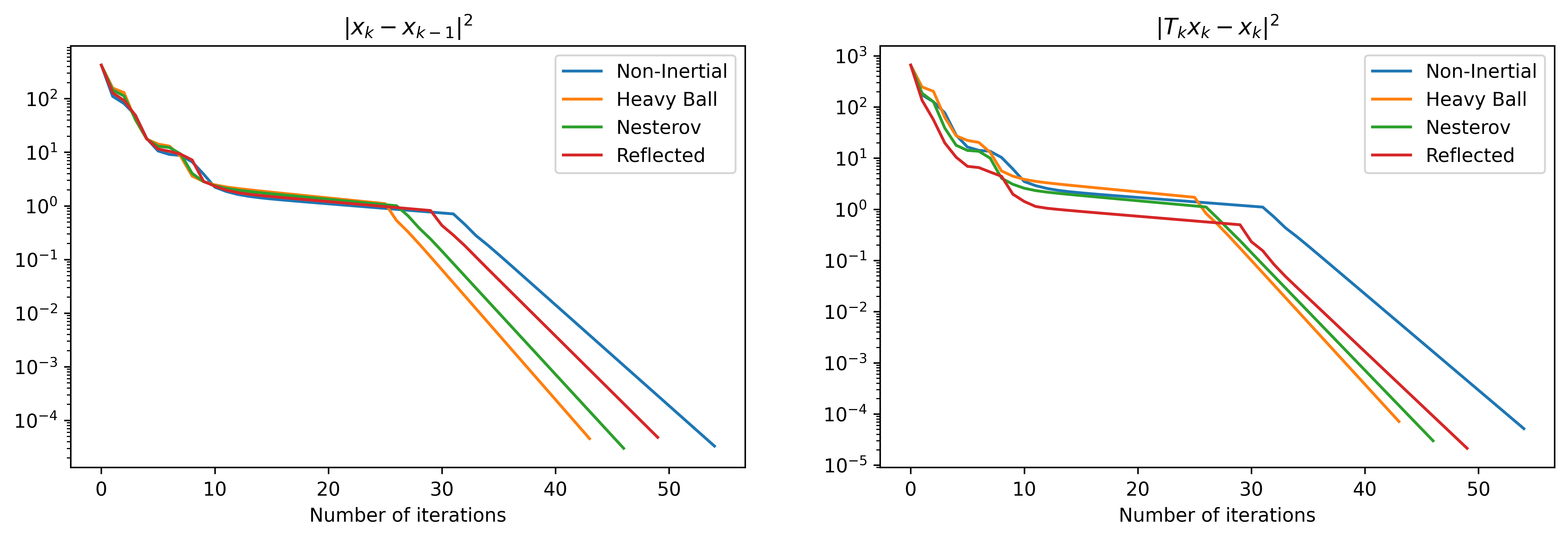}
    \caption{Convergence plots of the residuals for $\rho=2/ L$ and $\lambda=0.8$.}
    \label{fig:residuals}
\end{figure}

Figure \ref{fig:residuals} shows convergence rates of the residual quantities, using the parameters $\rho=\frac2L$ and $\lambda=0.8$, observed favorable throughout the previous two experiments. As expected, we do observe linear convergence asymptotically.

\section{Conclusions}

We have provided a systematic, unified and insightful analysis of the hypotheses that ensure the weak, strong and linear convergence of a family of inexact fixed point iterations combining relaxation with different inertial (acceleration) principles. Several previous results, obtained by analysing particular cases separately, are either contained or improved. The numerical illustrations reveal advantages of the use of inertia, and open questions about the optimality of the parameter choice.
{\small
\section{Acknowledgements}
The authors thank the anonymous reviewers for their constructive remarks, and in particular, for suggesting the tightness analysis in Subsection \ref{ssec:tight}. This research benefited from the support of the FMJH Program Gaspard Monge for optimization and operations research and their interactions with data science. }

{\small
\section{Data Availability}
Data sets generated during the current study are available from the corresponding author on reasonable request. }

\appendix
\section{Proofs of Proposition \ref{prop:2} and Theorem \ref{thm:strong}, Postponed from Section \ref{sec:convergence}}\label{sec:postponed}

The properties of the norm stated as Lemmas \ref{L:trick1} and \ref{L:trick2}, as well as those of real sequences given in Lemmas \ref{lemma:bounded} and \ref{L:diff}, all in Appendix \ref{sec:aux}, will be useful to shorten some of the proofs in the upcoming subsections.

\begin{remark} \label{remark:simplify}
Algorithm \eqref{alg:0}, as well as the convergence results stated in Theorems \ref{thm:weak} and \ref{thm:strong} , we may assume without loss of generality that $\theta_k\equiv 0$. Indeed, by setting $\tilde x_k=x_k-\theta_{k-1}$, (with $\theta_{-1}=\theta_0=0$), Algorithm \eqref{alg:0} is equivalent to
\begin{equation*}
\begin{cases}
y_k &= \quad \tilde x_k+\alpha_k(\tilde x_k-\tilde x_{k-1}) + \tilde \varepsilon_k,\qqbox{where}\tilde\varepsilon_k\coloneqq \theta_{k-1}+\alpha_k(\theta_{k-1}-\theta_{k-2})+ \varepsilon_k \\
z_k &= \quad \tilde x_k+\beta_k(\tilde x_k-\tilde x_{k-1}) + \tilde\rho_k,\qqbox{where}\tilde\rho_k\coloneqq\theta_{k-1}+\beta_k(\theta_{k-1}-\theta_{k-2})+ \rho_k \\
\tilde x_{k+1} &= \quad (1-\lambda_k)y_k+\lambda_kT_kz_k.
\end{cases}
\end{equation*}
Since $(\alpha_k)$ and $(\beta_k)$ are bounded, $(\tilde\varepsilon_k),(\tilde\rho_k)\in\ell^r(\H)$ whenever $(\varepsilon_k),(\rho_k), (\theta_k)\in\ell^r(\H)$.
\end{remark}

\subsection{Proof of Proposition \ref{prop:2}}\label{sec:prop2}

In order to simplify the notation, given $p\in \H$, we define, for $k\ge 1$,
\begin{equation} \label{E:Notation_1}
\begin{cases}
~~\nu_k &\coloneqq\quad \lambda_k^{-1}-1, \\
~~\mu_k &\coloneqq \quad (1-\lambda_k)\alpha_k+\lambda_k\beta_k, \\
~~\Delta_k(p) &\coloneqq \quad \|x_k-p\|^2-\|x_{k-1}-p\|^2, \\
~~B_k&\coloneqq \quad(1-\lambda_k)\alpha_k(1+\alpha_k) +\lambda_k\beta_k(1+\beta_k) +\nu_k\alpha_k(1-\alpha_k), \\
~~E_k^1 &\coloneqq \quad (1-\lambda_k - \nu_k)\|\varepsilon_k\|^2+\lambda_k\|\rho_k\|^2, \\
~~E_k^2&\coloneqq \quad (1-\lambda_k+\nu_k)\alpha_k\|\varepsilon_k\|+\lambda_k\beta_k\|\rho_k\|, \\
~~E_k^3&\coloneqq \quad (1-\lambda_k+\nu_k)(1+\alpha_k)\|\varepsilon_k\|+\lambda_k(1+\beta_k)\|\rho_k\|,\\
~~E_k^4&\coloneqq \quad \nu_k\|\varepsilon_k\|, \\
~~E_k(p) &\coloneqq \quad E_k^1+2E_k^2\|x_{k-1}-p\|+2E_k^3\|x_k-p\|+2E_k^4\|x_{k+1}-p\|.
\end{cases}
\end{equation}
We begin by proving the following:

\begin{lemma} \label{lemma:1}
Let $(\alpha_k), (\beta_k)$ be sequences in $[0,1]$,  $(\lambda_k)$ in $(0,1)$, and $(\varepsilon_k)$ and $ (\rho_k)$ in $\H$. Let $T_k\colon \H\to \H$ be a family of quasi-nonexpansive operators such that $F\coloneqq \bigcap_{k\ge 1}\Fix(T_k)\neq \emptyset$. Also, let $(x_k, y_k, z_k)$ be generated by Algorithm \eqref{alg:0}. Then, for all $p\in F$ and all $k\ge 1$, we have
\begin{equation*}
\Delta_{k+1}(p)\le
\mu_k\Delta_k(p)+B_k\|x_k-x_{k-1}\|^2 -\nu_k(1-\alpha_k)\|x_{k+1}-x_k\|^2 +E_k(p).
\end{equation*}
\end{lemma}

\begin{proof}
As explained in Remark \ref{remark:simplify}, we may assume $\theta_k\equiv 0$. Fix some $p\in F$. Using the definition of $x_k$, Lemma \ref{L:trick1} with $\zeta=-\lambda_k$, and the quasi-nonexpansiveness of $T_k$, it follows that
\begin{align}\label{eq:lemfirst}
\|x_{k+1}-p\|^2 &=\left\|(1-\lambda_k)y_k+\lambda_kT_kz_k-p\right\|^2 \nonumber \\
&= (1-\lambda_k)\left\|y_k-p\right\|^2-\lambda_k(1-\lambda_k)\left\|T_kz_k-y_k\right\|^2 + \lambda_k\|T_kz_k-p\|^2 \nonumber \\
&\le (1-\lambda_k)\left\|y_k-p\right\|^2-\lambda_k(1-\lambda_k)\left\|T_kz_k-y_k\right\|^2 + \lambda_k\|z_k-p\|^2.
\end{align}
Using the definition of $y_k$ and Lemma \ref{L:trick1} with $\zeta=\alpha_k$, we can bound the first term on the right-hand side by observing that
\begin{align}\label{eq:lemsec}
\left\|y_k-p\right\|^2 
&= \|x_k-p+\alpha_k(x_k-x_{k-1})+\varepsilon_k\|^2  \nonumber \\
&\le \|x_k-p+\alpha_k(x_k-x_{k-1})\|^2+2\|\varepsilon_k\| \|x_k-p+\alpha_k(x_k-x_{k-1})\|+\|\varepsilon_k\|^2  \nonumber \\
&= (1+\alpha_k)\|x_k-p\|^2+\alpha_k(1+\alpha_k)\|x_k-x_{k-1}\|^2-\alpha_k\|x_{k-1}-p\|^2 \nonumber \\
&\quad +2\|\varepsilon_k\| \|x_k-p+\alpha_k(x_k-p+p-x_{k-1})\|+\|\varepsilon_k\|^2  \nonumber \\
&\le (1+\alpha_k)\|x_k-p\|^2+\alpha_k(1+\alpha_k)\|x_k-x_{k-1}\|^2-\alpha_k\|x_{k-1}-p\|^2 \nonumber \\
&\quad +2\|\varepsilon_k\|\big[ (1+\alpha_k)\|x_k-p\|+\alpha_k\|x_{k-1}-p\|\big]+\|\varepsilon_k\|^2.
\end{align}
Analogously, the last term of \eqref{eq:lemfirst} may be bounded by
\begin{align}\label{eq:lemsecbis}
\|z_k-p\|^2 &\le ~(1+\beta_k)\|x_k-p\|^2+\beta_k(1+\beta_k)\|x_k-x_{k-1}\|^2-\beta_k\|x_{k-1}-p\|^2 \nonumber \\
&\quad +2\|\rho_k\|\big[(1+\beta_k)\|x_k-p\|+\beta_k\|x_{k-1}-p\|\big]+\|\rho_k\|^2.
\end{align}
For the middle term of \eqref{eq:lemfirst}, we use again the definition of $x_k$ and $y_k$ and Lemma \ref{L:trick1} with $\zeta=-\alpha_k$, to write
\begin{equation*}
\begin{split}
-\lambda_k^2\|T_kz_k-y_k\|^2
&=-\|x_{k+1}-y_k\|^2 \\
&=-\|x_{k+1}-x_k-\alpha_k(x_k-x_{k-1})-\varepsilon_k\|^2 \\
&= -\|x_{k+1}-x_{k}-\alpha_k(x_k-x_{k-1})\|^2 +2\|\varepsilon_k\|\|x_{k+1}-x_k-\alpha_k(x_k-x_{k-1})\|-\|\varepsilon_k\|^2 \\
&= -(1-\alpha_k)\|x_{k+1}-x_{k}\|^2+\alpha_k(1-\alpha_k)\|x_k-x_{k-1}\|^2 -\alpha_k\|x_{k+1}-2x_k+x_{k-1}\|^2\\
&\quad  +2 \|\varepsilon_k\|\|x_{k+1}-x_k-\alpha_k(x_k-x_{k-1})\|-\|\varepsilon_k\|^2. \\
\end{split}
\end{equation*}
We disregard the third-to-last term on the right-hand side, and use the same trick as above on the coefficient of $\|\varepsilon_k\|$, to deduce that
\begin{align}\label{eq:lemthi}
-\lambda_k^2\|y_k-T_kz_k\|^2&\le-(1-\alpha_k)\|x_{k+1}-x_{k}\|^2 +\alpha_k(1-\alpha_k)\|x_k-x_{k-1}\|^2 \nonumber \\
&\quad +2\|\varepsilon_k\|\big[\|x_{k+1}-p\|+(1+\alpha_k)\|x_k-p\|+\alpha_k\|x_{k-1}-p\|\big] - \|\varepsilon_k\|^2.
\end{align} 
Combining inequalities \eqref{eq:lemfirst}, \eqref{eq:lemsec}, \eqref{eq:lemsecbis} and \eqref{eq:lemthi}, 
we obtain
\begin{align*}
\|x_{k+1}-p\|^2 \le &\ (1-\lambda_k)\left\|y_k-p\right\|^2-\nu_k\lambda_k^2\left\|T_kz_k-y_k\right\|^2 + \lambda_k\|z_k-p\|^2  \\
\le &\ (1-\lambda_k)\bigg[(1+\alpha_k)\|x_k-p\|^2+\alpha_k(1+\alpha_k)\|x_k-x_{k-1}\|^2 -\alpha_k\|x_{k-1}-p\|^2\nonumber \\
&\qquad +2\|\varepsilon_k\|\big[ (1+\alpha_k)\|x_k-p\|+\alpha_k\|x_{k-1}-p\|\big]+\|\varepsilon_k\|^2 \bigg] \\
&\  +\nu_k\bigg[-(1-\alpha_k)\|x_{k+1}-x_{k}\|^2 +\alpha_k(1-\alpha_k)\|x_k-x_{k-1}\|^2 \nonumber \\
&\qquad +2\|\varepsilon_k\|\big[\|x_{k+1}-p\|+(1+\alpha_k)\|x_k-p\|+\alpha_k\|x_{k-1}-p\|\big] - \|\varepsilon_k\|^2\bigg] \\
&\  +\lambda_k\bigg[(1+\beta_k)\|x_k-p\|^2+\beta_k(1+\beta_k)\|x_k-x_{k-1}\|^2 -\beta_k\|x_{k-1}-p\|^2 \\
&\qquad +2\|\rho_k\|\big[(1+\beta_k)\|x_k-p\|+\beta_k\|x_{k-1}-p\|\big]+\|\rho_k\|^2\bigg].
\end{align*} 

Subtracting $\|x_k-p\|^2$ on both sides and aggregating similar terms, we get
\begin{align*}
\Delta_{k+1}(p) 
\le &\ \bigg[(1-\lambda_k)(1+\alpha_k)+\lambda_k(1+\beta_k)-1\bigg]\|x_k-p\|^2-\bigg[(1-\lambda_k)\alpha_k+\lambda_k\beta_k\bigg]\|x_{k-1}-p\|^2  \nonumber \\
&\  +\bigg[(1-\lambda_k)\alpha_k(1+\alpha_k)+\lambda_k\beta_k(1+\beta_k) +\nu_k\alpha_k(1-\alpha_k)\bigg]\|x_k-x_{k-1}\|^2 -\nu_k(1-\alpha_k)\|x_{k+1}-x_{k}\|^2\\
&\ +\bigg[(1-\lambda_k-\nu_k)\|\varepsilon_k\|^2+\lambda_k\|\rho_k\|^2\bigg]
+2\bigg[(1-\lambda_k+\nu_k)\alpha_k\|\varepsilon_k\|+\lambda_k\beta_k\|\rho_k\|\bigg]\|x_{k-1}-p\| \\
&\ +2\bigg[(1-\lambda_k+\nu_k)(1+\alpha_k)\|\varepsilon_k\|+\lambda_k(1+\beta_k)\|\rho_k\|\bigg]\|x_k-p\| +2\nu_k\|\varepsilon_k\|\|x_{k+1}-p\|,
\end{align*} 
which completes the proof.
\qed
\end{proof}

We may now return to the proof of Proposition \ref{prop:2}.

As before (by Remark \ref{remark:simplify}), we assume $\theta_k\equiv 0$. Fix $p\in F$.
Inequality \eqref{eq:tbr} implies that there is $\rho>0$ such that
\begin{equation*}
B_k\le \nu_{k-1}(1-\alpha_{k-1})-\rho
\end{equation*}
for all $k\ge 1$.
Combining this with Lemma \ref{lemma:1}, and writing $\delta_k\coloneqq \nu_{k-1}(1-\alpha_{k-1})\|x_k-x_{k-1}\|^2$, we obtain
\begin{equation} \label{eq:long}
\Delta_{k+1}(p)
\le \mu_k\Delta_k(p) +\delta_k -\delta_{k+1} -\rho\|x_k-x_{k-1}\|^2  +E_k(p).
\end{equation}
Summing for $k=1,\ldots, j$, and recalling the definitions in \eqref{E:Notation_1}, and the fact that $\mu_k$ is nondecreasing and bounded above by $M<1$, we obtain
$$\|x_{j+1}-p\|^2\le M\|x_j-p\|^2+\left[A+\sum_{k=1}^\infty E_k^1\right]+\sum_{k=1}^j \left[\frac{}{}\!E_k^2\|x_{k-1}-p\|+E_k^3\|x_k-p\|+E_k^4\|x_{k+1}-p\|\right],$$
where $A$ collects the constant terms from the telescopic sum. In other words, 
$$\zeta_{j+1}^2-M\zeta_j^2\le C+\sum_{k=1}^j\sum_{i=0}^{k+1}e_{k,i}\zeta_i,$$
where $\zeta_j=\|x_j-p\|$, $C=A+\sum_{k=1}^\infty E_k^1$, $e_{k,i}=E_k^{i+3-k}$ for $i=k-1,k,k+1$ and $e_{k,i}=0$ if $i<k-1$. Lemma \ref{lemma:bounded} shows that $\|x_k-p\|$ is bounded, and so $E_k(p)$ is summable. On the other hand, \eqref{eq:long} also implies
$$\rho\|x_k-x_{k-1}\|^2\le \left[\|x_k-p\|^2-\mu_{k-1}\|x_{k-1}-p\|^2+\delta_k\right]- \left[\|x_{k+1}-p\|^2-\mu_k\|x_k-p\|^2+\delta_{k+1}\right]+ E_k(p).$$
Since $E_k(p)$ is summable and the sums in the brackets are bounded from below, 
$\sum_{k=1}^\infty \|x_k-x_{k-1}\|^2$ is convergent. Considering that
\begin{equation*}
\lambda_k^2\|T_kz_k-y_k\|^2\le 4\|x_{k+1}-x_{k}\|^2 +4 \alpha_k^2\|x_k-x_{k-1}\|^2+2\|\varepsilon_k\|^2,
\end{equation*}
$\sum_{k=1}^\infty \|T_kz_k-y_k\|^2$ must converge as well. 

Finally, to prove that $\|x_k-p\|$ is convergent, note that \eqref{eq:long} implies that $\Delta_{k+1}(p)\le \mu_k\Delta_{k}(p)+\delta_k+E_k(p)$. It therefore suffices to apply Lemma \ref{L:diff} with $\Omega_k=\Delta_k(p)$, $\omega_k=\|x_k-p\|^2$, $b=1$, $a_k=\mu_k$ and $d_k=\delta_k+E_k(p)$ to conclude. \qed

\subsection{Proof of Theorem \ref{thm:strong}}\label{sec:thmstrong}

As explained in Remark \ref{remark:simplify}, we may assume $\theta_k\equiv 0$. As in the previous section, we simplify the notation by setting
\begin{equation} \label{E:Notation_2}
\begin{cases}
~~\nu_k &\coloneqq \quad \lambda_k^{-1}-1, \\
~~\mu_k &\coloneqq \quad (1-\lambda_k)\alpha_k+\lambda_kq_k^2\beta_k, \\
~~Q_k &\coloneqq \quad 1-\lambda_k+\lambda_kq_k^2, \\
~~\mathcal{B}_k(\gamma_1, \gamma_2)&\coloneqq \quad (1+\gamma_1)\biggl[(1-\lambda_k)\alpha_k(1+\alpha_k) +\lambda_kq_k^2\beta_k(1+\beta_k)\biggr] +(1-\gamma_2)\nu_k\alpha_k(1-\alpha_k), \\
~~\mathcal E_k(\gamma_1, \gamma_2) &\coloneqq \quad (\gamma_1^{-1}+1)\biggl[(1-\lambda_k)\|\varepsilon_k\|^2+\lambda_kq_k^2\|\rho_k\|^2\biggr] +\nu_k(\gamma_2^{-1}-1)\|\varepsilon_k\|^2.
\end{cases}
\end{equation}
Since the definition of $\mu_k$ in \eqref{E:Notation_1} corresponds exactly to the case where $q_k\equiv 1$, using the same notation should not lead to confusion. If $\varepsilon_k\equiv\rho_k\equiv 0$, we allow $\gamma_1=\gamma_2=0$ and define $\mathcal E_k(0, 0)=0$.

By Lemma \ref{L:trick1} with $u=y_k-p^*$, $v=T_kz_k-y_k$ and $\zeta=-\lambda_k$, and the $q_k$-quasi-contractivity of $T_k$, we have
\begin{equation}\label{prop:conv1}
\begin{split}
\|x_{k+1}-p^*\|^2 
&=\|y_k-p^*+\lambda_k(T_kz_k-y_k)\|^2 \\
&=(1-\lambda_k)\|y_k-p^*\|^2+\lambda_k(\lambda_k-1)\|T_kz_k-y_k\|^2+\lambda_k\|T_kz_k-p^*\|^2 \\
&\le (1-\lambda_k)\|y_k-p^*\|^2+\lambda_k(\lambda_k-1)\|T_kz_k-y_k\|^2+\lambda_kq_k^2\|z_k-p^*\|^2.
\end{split}
\end{equation}
Applying Lemma \ref{L:trick2} with $\gamma=\gamma_1>0$ and Lemma \ref{L:trick1} with $\zeta=\alpha_k$, we obtain
\begin{align}\label{eq2:lemsec}
\left\|y_k-p^*\right\|^2 
&= \|x_k-p^*+\alpha_k(x_k-x_{k-1})+\varepsilon_k\|^2  \nonumber \\
&\le (1+\gamma_1)\|x_k-p^*+\alpha_k(x_k-x_{k-1})\|^2+(1+\gamma_1^{-1})\|\varepsilon_k\|^2  \nonumber \\
&= (1+\gamma_1)\left[(1+\alpha_k)\|x_k-p^*\|^2+\alpha_k(1+\alpha_k)\|x_k-x_{k-1}\|^2-\alpha_k\|x_{k-1}-p^*\|^2\right] +(1+\gamma_1^{-1})\|\varepsilon_k\|^2.
\end{align}
Analogously, 
\begin{align}\label{eq2:lemsecbis}
\left\|z_k-p^*\right\|^2 
&\le (1+\gamma_1)\left[(1+\beta_k)\|x_k-p^*\|^2+\beta_k(1+\beta_k)\|x_k-x_{k-1}\|^2-\beta_k\|x_{k-1}-p^*\|^2\right] +(1+\gamma_1^{-1})\|\rho_k\|^2.
\end{align}
For the middle term of \eqref{prop:conv1}, we use once again the definition of $x_k$ and $y_k$, Lemma \ref{L:trick2} with $\gamma=\gamma_2>0$ and Lemma \ref{L:trick1} with $\zeta=-\alpha_k$, to write
\begin{equation*}
\begin{split}
-\lambda_k^2\|T_kz_k-y_k\|^2
&=-\|x_{k+1}-y_k\|^2 \\
&=-\|x_{k+1}-x_k-\alpha_k(x_k-x_{k-1})-\varepsilon_k\|^2 \\
&= -(1-\gamma_2)\|x_{k+1}-x_{k}-\alpha_k(x_k-x_{k-1})\|^2-(1-\gamma_2^{-1})\|\varepsilon_k\|^2 \\
&= (1-\gamma_2)\left[(1-\alpha_k)\|x_{k+1}-x_{k}\|^2+\alpha_k(1-\alpha_k)\|x_k-x_{k-1}\|^2 -\alpha_k\|x_{k+1}-2x_k+x_{k-1}\|^2\right]\\
&\qquad  -(1-\gamma_2^{-1})\|\varepsilon_k\|^2. \\
\end{split}
\end{equation*}
We multiply this equation by $\nu_k$, and disregard the second-to-last term on the right-hand side, to rewrite it as
\begin{align}\label{eq2:lemthi}
-\lambda_k(1-\lambda_k)\|y_k-T_kz_k\|^2&\le-(1-\gamma_2)\nu_k(1-\alpha_k)\|x_{k+1}-x_{k}\|^2 +(1-\gamma_2)\nu_k\alpha_k(1-\alpha_k)\|x_k-x_{k-1}\|^2 \nonumber\\
&\qquad + (1-\gamma_2^{-1})\nu_k\|\varepsilon_k\|^2.
\end{align} 
We combine Inequalities \eqref{eq2:lemsec}, \eqref{eq2:lemsecbis} and \eqref{eq2:lemthi} with \eqref{prop:conv1}, to deduce that
\begin{equation*}
\begin{split}
\|x_{k+1}-p^*\|^2  \le  ~
&(1+\gamma_1) \biggl[(1-\lambda_k)(1+\alpha_k)+\lambda_kq_k^2(1+\beta_k)\biggr]\|x_k-p^*\|^2 \\
&\quad -(1+\gamma_1)\biggl[(1-\lambda_k)\alpha_k+\lambda_kq_k^2\beta_k\biggr]\|x_{k-1}-p^*\|^2 +\mathcal B_k(\gamma_1, \gamma_2)\|x_k-x_{k-1}\|^2  \\
&\quad  -(1-\gamma_2)\nu_k(1-\alpha_k)\|x_{k+1}-x_k\|^2  +\mathcal E_k(\gamma_1, \gamma_2),
\end{split}
\end{equation*}
for all $\gamma_1,\gamma_2>0$. Using the definitions of $\mu_k$ and $Q_k$, we rewrite this as
\begin{equation} \label{eq:anchor_strong}
\begin{split}
\|x_{k+1}-p^*\|^2  \le ~
&(1+\gamma_1)(\mu_k+Q_k)\|x_k-p^*\|^2 -(1+\gamma_1)\mu_k\|x_{k-1}-p^*\|^2 +\mathcal B_k(\gamma_1, \gamma_2)\|x_k-x_{k-1}\|^2\\
&\quad -(1-\gamma_2)\nu_k(1-\alpha_k)\|x_{k+1}-x_k\|^2 +\mathcal E_k(\gamma_1, \gamma_2).
\end{split}
\end{equation}
Since $Q=\sup Q_k\le 1-\lambda(1-q^2)<1$, we can select $\gamma_1>0$ such that $(1+\gamma_1)Q<1$. Since Inequality \eqref{eq:tbr_strong} remains valid if we multiply it by $(1+\gamma_1)$, we have
$$\sup_{k\ge 1}\left[\frac{}{}\!\mathcal B_k(\gamma_1, \gamma_2)-(1+\gamma_1)Q_k\nu_{k-1}(1-\alpha_{k-1})\right]<0,$$
so we can pick $\gamma_2\in (0,1)$ such that 
\begin{equation}\label{eq:ass}
\mathcal B_k(\gamma_1, \gamma_2) \le (1+\gamma_1)(1-\gamma_2)Q_k\nu_{k-1}(1-\alpha_{k-1})
\end{equation}
for all $k\ge 1$. Using \eqref{eq:ass} in \eqref{eq:anchor_strong}, we get
\begin{equation*}
\begin{split}
\|x_{k+1}-p^*\|^2  \le ~
&(1+\gamma_1)(\mu_k+Q_k)\|x_k-p^*\|^2 -(1+\gamma_1)\mu_k\|x_{k-1}-p^*\|^2  +\mathcal E_k(\gamma_1, \gamma_2)  \\
&\quad +(1+\gamma_1)(1-\gamma_2)Q_k\nu_{k-1}(1-\alpha_{k-1})\|x_k-x_{k-1}\|^2  -(1-\gamma_2)\nu_k(1-\alpha_k)\|x_{k+1}-x_k\|^2.
\end{split}
\end{equation*}
Since $(1+\gamma_1)Q<1$ and $(\mu_k)$ is nondecreasing, we can group some terms to deduce that
\begin{equation*}
\begin{split}
& \|x_{k+1}-p^*\|^2 -(1+\gamma_1)\mu_k\|x_k-p^*\|^2 + (1-\gamma_2)\nu_k(1-\alpha_k)\|x_{k+1}-x_k\|^2 \\
& \quad \le 
 (1+\gamma_1)Q_k\left[\frac{}{}\!\|x_k-p^*\|^2 -(1+\gamma_1)\mu_{k-1}\|x_{k-1}-p^*\|^2 +(1-\gamma_2)\nu_{k-1}(1-\alpha_{k-1})\|x_k-x_{k-1}\|^2\right] +\mathcal E_k(\gamma_1, \gamma_2).
\end{split}
\end{equation*}
Defining 
\begin{equation*}
C_k\coloneqq \|x_{k}-p^*\|^2 - (1+\gamma_1)\mu_{k-1}\|x_{k-1}-p^*\|^2 + (1-\gamma_2)\nu_{k-1}(1-\alpha_{k-1})\|x_{k}-x_{k-1}\|^2,
\end{equation*}
this reads
\begin{equation}\label{eq:random}
C_{k+1}  \le (1+\gamma_1)Q_k C_k  +\mathcal E_k(\gamma_1, \gamma_2).
\end{equation}

We can now use Lemma \ref{L:diff} with $\Omega_k=C_k$, $\omega_k=\|x_{k}-p^*\|^2$, $b=(1+\gamma_1)M<1$, $a_k=(1+\gamma_1)Q_k$, $a=(1+\gamma_1)Q<1$ and $d_k=\mathcal E_k(\gamma_1,\gamma_2)$, to deduce that $\sum \|x_k-p^*\|^2<+\infty$, and so $\|x_k-p^*\|^2$ converges to $0$. This implies the wanted convergence of $(x_k)$ to $p^*$, and hence, since $(\alpha_k)$ and $(\beta_k)$ are bounded and $\varepsilon_k$ and $\rho_k$ converge to $0$, we conclude the convergence of $(x_k, y_k, z_k)$ to $(p^*, p^*, p^*)$. Now set $D_k=\min_{j=1,\dots,k}\|x_j-p^*\|^2$, for $k\ge 1$. We have
$$kD_k\le\sum_{j=1}^kD_j\le\sum_{j=1}^k\|x_j-p^*\|^2\le\frac{1}{\big(1-(1+\gamma_1)M\big)\big(1-(1+\gamma_1)Q\big)}\left[2\|x_1-p^*\|^2+\sum_{k=1}^\infty \mathcal E_k(\gamma_1,\gamma_2)\right],$$  
and so $D_k=\mathcal O\big(\frac{1}{k}\big)$. Moreover, since the sequence $(D_k)$ is nonincreasing and belongs to $\ell^1(\R)$, we have $\limty{k}kD_k=0$. In other words, $D_k=o\left(\frac{1}{k}\right)$ as $k\to\infty$.

In the absence of perturbations ($\varepsilon_k\equiv \rho_k\equiv 0$), the sequence $x_k$ converges linearly to $p^*$. Indeed, in that case, $\mathcal E_k(\gamma_1,\gamma_2)\equiv 0$, and we can take $\gamma_1=0$, so that \eqref{eq:random} reduces to $C_{k+1}\le Q_kC_k$. Using Lemma \ref{L:diff} with $\Omega_k=C_k$, $a_k=Q_k$, $a=Q<1$ and $d_k\equiv 0$, we see that $[C_{k+1}]_+\le Q^k[C_1]_+=Q^k\|x_1-p^*\|^2$. 
Using the definition of $C_k$ and iterating, one obtains 
 \begin{equation*}
     \|x_{k+1}-p^*\|^2\le M\|x_k-p^*\|^2+Q^{k}\|x_1-p^*\|^2\le \cdots \le \left[\sum_{i=0}^{k}M^iQ^{k-i}\right]\|x_1-p^*\|^2= \frac{Q^{k+1}-M^{k+1}}{Q-M}\|x_1-p^*\|^2.
\end{equation*}
From \eqref{E:Notation_2}, we see that $\mu_k\le Q_k+(1-\Lambda)(A-1)$ for all $k\ge 1$, whence $Q-M\ge(1-\Lambda)(1-A)>0$, and the conclusion follows.
\qed 
\section{Auxiliary Results}\label{sec:aux}

The following properties of the norm were used repeatedly in Appendix \ref{sec:postponed}. The first one is a generalized parallelogram identity:

\begin{lemma} \label{L:trick1}
For every $u,v\in\H$ and $\zeta\in\R$, we have
\begin{equation*}
\|u+\zeta v\|^2 = (1+\zeta)\|u\|^2+\zeta(1+\zeta)\|v\|^2-\zeta\|u-v\|^2.
\end{equation*}
\end{lemma}

\begin{proof}
It suffices to add the identities
\begin{eqnarray*}
\|u+\zeta v\|^2 & = & \|u\|^2+\zeta^2\|v\|^2+2\zeta\langle u,v\rangle \\  
\zeta\|u-v\|^2 & = & \zeta\|u\|^2+\zeta\|v\|^2-2\zeta\langle u,v\rangle
\end{eqnarray*}
and rearrange the terms.
\end{proof}

The next one is a direct consequence of Cauchy-Schwarz and Young's inequalities:

\begin{lemma} \label{L:trick2}
For every $u,v\in\H$ and $\gamma>0$, we have
\begin{equation*}
 \left(1-\gamma\right)\|u\|^2+\left(1-\frac1\gamma\right)\|v\|^2\le \|u\pm v\|^2\le \left(1+\gamma\right)\|u\|^2+\left(1+\frac1\gamma\right)\|v\|^2.
\end{equation*}    
\end{lemma}

\begin{proof}
We first bound
$$\left|\frac{}{}\!\|u \pm v\|^2-\|u\|^2-\|v\|^2\right|=2|\langle u,v\rangle|\le 2\left(\frac{}{}\!\sqrt{\gamma}\|u\|\right)\left(\frac{1}{\sqrt{\gamma}}\|v\|\right)\le \gamma\|u\|^2+\frac{1}{\gamma}\|v\|^2,$$
and then rewrite without the absolute value.
\end{proof}

We now provide two elementary but not so standard results on real sequences that have been used in Section \ref{sec:convergence}. The first one is an extension of \cite[Lemma 5.14]{attouch2018fast}. 

\begin{lemma}\label{lemma:bounded}
Let $(\zeta_k)$ be a nonnegative sequence such that 
\begin{equation} \label{E:lemma_bounded}
\zeta_{j+1}^2-M\zeta_j^2\le C+\sum_{k=1}^j\sum_{i=0}^{k+1}e_{k,i}\zeta_i,    
\end{equation}
where $C\ge 0$, $0\le M<1$ and $e_{k,i}\ge 0$, with $D:=\sum_{k=1}^\infty\sum_{i=0}^{k+1} e_{k,i}<\infty$ for each $i$. Then, $(\zeta_j)$ is bounded.
\end{lemma}

\begin{proof}
Set $Z_j=\max_{i=1,\dots j}\zeta_i$. Take $n\ge 1$. For every $j\le n$, \eqref{E:lemma_bounded} gives 
$$\zeta_{j+1}^2\le MZ_{n}^2+C+Z_{n+1}\sum_{k=1}^j\sum_{i=0}^{k+1}e_{k,i}\le  MZ_{n+1}^2+C+DZ_{n+1}.$$
Since the right-hand side does not depend on $j$, we may take the maximum for $j=1,\dots,n$, and rearrange the terms, to obtain
$$(1-M)Z_{n+1}^2-DZ_{n+1}-C\le 0,$$
which implies $(\zeta_j)$ is bounded.
\end{proof}

Now, write $[\omega]_+=\max\{0, \omega\}$ for $\omega\in\R$. The following result is a straightforward extension of \cite[Lemma 2.2]{mainge2008convergence}, a result that was actually established much earlier, embedded in the proof of \cite[Theorem 2.1]{alvarez2001inertial}.

\begin{lemma}\label{L:diff}
    Let $(a_k)$ and $(d_k)$ be nonnegative sequences such that $a_k\le a<1$ for all $k\ge 0$, and $\sum_{k=0}^\infty d_k<+\infty$. Consider a real sequence $(\Omega_k)$ such that
    $$\Omega_{k+1}\le a_k\Omega_k+d_k$$
    for all $k\ge 0$. Then $\sum_{k=1}^\infty [\Omega_k]_{+}$ is convergent. If, moreover, $d_k\equiv 0$, then $[\Omega_{k+1}]_+\le a^k[\Omega_1]_+$. Either way, if $\Omega_k\ge \omega_k-b\omega_{k-1}$ for all $k\ge 0$, where $b\in [0, 1]$ and $(\omega_k)$ is nonnegative, then $(\omega_k)$ is convergent. If, moreover, $b<1$, then $\sum_{k=0}^\infty \omega_k<+\infty$.
\end{lemma}

\bibliographystyle{abbrv}
\bibliography{ref.bib}

\end{document}